\documentclass{fic-l}
\usepackage{epsf}
\usepackage[dvips]{graphics} 
\usepackage{amssymb,amsfonts,latexsym}
\newtheorem{theorem}{Theorem}[section] 
\newtheorem{lemma}[theorem]{Lemma}
\newtheorem{proposition}[theorem]{Proposition}
\theoremstyle{definition}
\newtheorem{definition}[theorem]{Definition}

\theoremstyle{remark} \newtheorem{remark}[theorem]{Remark}
\numberwithin{equation}{section}
\catcode`\@=11
\newif\if@fewtab\@fewtabtrue
{\count255=\time\divide\count255 by 60
\xdef\hourmin{\number\count255}
\multiply\count255 by-60\advance\count255 by\time
\xdef\hourmin{\hourmin:\ifnum\count255<10 0\fi\the\count255}}
\def\ps@draft{\let\@mkboth\@gobbletwo
    \def\@oddfoot{\hbox to 7 cm{\tiny \versionno
       \hfil}\hskip -7cm\hfil\rm\thepage \hfil {\tiny\draftdate}}
    \def\@oddhead{}
    \def\@evenhead{}\let\@evenfoot\@oddfoot}
\def\draftdate{\number\month/\number\day/\number\year\ \ \ \hourmin }

\def\citen#1{\if@filesw \immediate\write \@auxout {\string\citation{#1}}\fi%
\@tempcntb\m@ne \let\@h@ld\relax \def\@citea{}%
\@for \@citeb:=#1\do {\@ifundefined {b@\@citeb}%
    {\@h@ld\@citea\@tempcntb\m@ne{\bf ?}%
    \@warning {Citation `\@citeb ' on page \thepage \space undefined}}%
    {\@tempcnta\@tempcntb \advance\@tempcnta\@ne
    \setbox\z@\hbox\bgroup\ifcat0\csname b@\@citeb \endcsname \relax
    \egroup \@tempcntb\number\csname b@\@citeb \endcsname \relax
    \else \egroup \@tempcntb\m@ne \fi \ifnum\@tempcnta=\@tempcntb
    \ifx\@h@ld\relax \edef \@h@ld{\@citea\csname b@\@citeb\endcsname}%
    \else \edef\@h@ld{\hbox{--}\penalty\@highpenalty
    \csname b@\@citeb\endcsname}\fi
    \else \@h@ld\@citea\csname b@\@citeb \endcsname \let\@h@ld\relax \fi}%
\def\@citea{,\penalty\@highpenalty\hskip.13em plus.13em minus.13em}}\@h@ld}
\def\@citex[#1]#2{\@cite{\citen{#2}}{#1}}%
\def\@cite#1#2{\leavevmode\unskip\ifnum\lastpenalty=\z@\penalty\@highpenalty\fi%
  \ [{\multiply\@highpenalty 3 #1%
  \if@tempswa,\penalty\@highpenalty\ #2\fi}]}   %
\makeatother 
\catcode`\@=12

\def\A             {\mbox{$A$}}
\def\Av            {{A\Vee}}
\def\abar          {\mbox{$\mathfrak A$}} 
\def\abulk         {\mbox{${\mathfrak A}_{\rm bulk}$}}
\def\absi          {absolutely simple}

\def\aext          {\mbox{${\mathfrak A}_{\rm ext}$}}
\def\alg           {algebra}
\def\auto          {automorphism}
\def\bc            {boundary condition}
\def\Bc            {Boundary condition}
\def\be            {\begin{equation}}
\def\bearl         {\begin{array}{l}}
\def\bearll        {\begin{array}{ll}}
\def\bearlll       {\begin{array}{lll}}
\def\bet           {\beta_{\!A}}

\def\calc          {\mbox{$\mathcal C$}}
\def\calca         {\mbox{${\mathcal C}_{\!A}$}}
\def\calcao        {\mbox{${\mathcal C}_{\!A}^0$}}
\def\calcm         {\mbox{${\mathcal C}_{{\rm Rep}\,\mathfrak A}$}}
\def\calco         {\mbox{${\mathcal C}^\circ$}}
\def\cald          {d} 
\def\cft           {conformal field theory}
\def\Cft           {Conformal field theory}
\def\cfts          {conformal field theories}
\def\Chi           {\mathcal X}
\def\chii          {\raisebox{.15em}{$\chi$}}
\def\cir           {\,{\circ}\,}
\def\class         {classification}
\def\coker         {{\rm Coker}\,}
\def\complex       {{\mathbb C}}
\def\con           {conformal }
\def\corfu         {correlation function}
\def\dd            {$\!\!$: }

\def\dim           {{\rm dim}}
\def\Dim           {{\rm Dim}}
\def\diml          {\dim_{\rm L}}
\def\dimL          {\dim}
\def\Diml          {\Dim_{\rm L}}
\def\dimr          {\dim_{\rm R}}
\def\dimR          {\dim}
\def\Dimr          {\Dim_{\rm R}}
\def\Dot           {\dot}
\def\dsty          {\displaystyle}
\newcommand\e[1]   {e_{#1}}

\def\eee           {\end{equation}}
\def\eear          {\end{array}}
\newcommand\eev[1] {{{}^\vee}\!{#1}}
\newcommand\Eev[1] {{{}^{\vee\!}}\!{#1}}
\newcommand\eew[1] {{{}^\wedge}\!{#1}}
\def\ej            {\e\J}
\def\ek            {\e\K}
\def\eps           {\epsilon}
\def\eq            {\,{=}\,}
\newcommand\erf[2] {(\ref{#1#2})}
\def\F             {{\sf F}}
\def\findim        {fi\-ni\-te-di\-men\-si\-o\-nal}

\def\G             {{\sf G}}
\def\gam           {\beta_I}
\def\hapl          {haploid}
\def\hapy          {haploidity}
\newcommand\hilde[1] {{{#1}}^\sharp}
\def\Hom           {{\rm Hom}}
\newcommand\hsp[1] {\mbox{\hspace{#1 em}}}
\def\hy            {$\mbox{-\hspace{-.66 mm}-}$}
\def\id            {\mbox{\sl id}}

\def\iN            {\,{\in}\,}
\def\inda          {{\rm Ind}_A}
\def\indar         {\widetilde{{\rm Ind}}_A}
\def\J             {{\rm J}}
\def\K             {{\rm K}}

\newcommand\labl[2]{\label{#1#2}}
\def\lie           {Lie algebra}
\def\llb           {\mbox{\large[}}
\def\Llb           {\mbox{\Large(}}
\def\lrb           {\mbox{\large]}}
\def\Lrb           {\mbox{\Large)}}
\def\mtc           {modular tensor category}
\def\Modinv        {Modular invarian}
\def\modinv        {modular invarian}
\def\Nu            {{\mathcal V}}
\def\obj           {{\mathcal O}bj}
\def\one           {I}
\def\onematrix     {\mbox{\small $1\!\!$}1}

\def\oa            {o} 
\def\ob            {{\check\oa}}

\def\ota           {\,{\otimes}_{\!A}^{}\,}
\def\Ota           {{\otimes}_{\!A}^{}}
\def\oti           {\,{\otimes}\,}
\def\Oti           {{\otimes}}
\def\parfu         {partition function}
\def\q             {quantum }
\def\Q             {Quantum }
\def\qft           {quantum field theory}
\def\qfts          {quantum field theories}
\def\R             {{\rho}}
\def\rep           {representation}
\def\Rep           {Representation}
\def\resp          {respectively}
\def\sa            {s}  
\def\sb            {{\check\sa}}
\def\smallone      {I}
\def\swap          {{sw\aa p}}
\def\swap          {{swap}}
\def\Swap          {{Swap}}
\def\swapc         {\swap-com\-mu\-ta\-ti\-ve}
\def\swapy         {\swap-com\-mu\-ta\-ti\-vi\-ty}
\def\Swapy         {\Swap-com\-mu\-ta\-ti\-vi\-ty}
\def\tc            {tensor category}
\def\tcc           {$\tilde c$-com\-mu\-ta\-ti\-ve}
\def\tccc          {$\tilde c$-cocom\-mu\-ta\-ti\-ve}
\def\tcy           {$\tilde c$-com\-mu\-ta\-ti\-vity}
\newcommand\tf[1]  {\tilde E(#1)}
\def\tft           {topological field theory}
\newcommand\tg[1]  {\tilde R(#1)}

\def\trl           {{\rm tr}_{\rm L}}
\def\Trl           {{\rm Tr}_{\rm L}}
\def\trr           {{\rm tr}_{\rm R}}
\def\Trr           {{\rm Tr}_{\rm R}}
\def\twodim        {two-di\-men\-si\-o\-nal}
\def\va            {v}  
\def\vb            {{\check\va}}

\def\Vee           {^\vee}

\def\voa           {vertex operator algebra}
\def\vsp           {{}\mbox{$\ $}\\[-1.8em]\phantom.}
\def\Wee           {^\wedge}
\newcommand\wilde[1] {{{#1}}^A}

\def\wrtt          {with respect to the }
\def\zet           {{\mathbb Z}}

\def\currentvolume{39}
\def\currentyear{25\,--\,71~(2003)}
\def\copyrightyear{2003}
\begin{document}
\setcounter{page}{25}
\title{Category theory for conformal boundary conditions}
\author{J\"urgen Fuchs}
\address{Institutionen f\"or fysik, Universitetsgatan 1,
S\,--\,651\,88\, Karlstad} \email{jfuchs@fuchs.tekn.kau.se}
\author{Christoph Schweigert}
\address{Fachbereich Mathematik, Schwerpunkt Algebra und Zahlentheorie\\
Universit\"at Hamburg, Bundesstra\ss e 55, D\,--\,20\,146 Hamburg}
\email{schweigert@math.uni-hamburg.de}
\date{}

\begin{abstract}
We study properties of the category of modules of an algebra object
$A$ in a tensor category $\mathcal C$. We show that the module category 
inherits various structures from $\mathcal C$, provided that $A$ is a 
Frobenius algebra with certain additional properties. As a by-product we 
obtain results about the Frobenius\hy Schur indicator in sovereign tensor 
categories. A braiding on $\mathcal C$ is not needed, nor is semisimplicity.
\\
We apply our results to the description of boundary conditions in 
two-dimensional conformal field theory and present illustrative examples. 
We show that when the module category is tensor, then it gives rise to a 
NIM-rep of the fusion rules, and discuss a possible relation with the 
representation theory of vertex operator algebras.
\end{abstract}

\maketitle

\section{CFT boundary conditions}

Boundary conditions in conformal field theory have various physical
applications, ranging from the study of defects in condensed matter 
physics to the theory of open strings. Such boundary conditions are 
partially characterized by the maximal vertex operator sub\alg\ \abar\ of 
the bulk chiral \alg\ \abulk\ that they respect \cite{fuSc112,Scfw}. That 
\abar\ is respected by a \bc\ means that the \corfu s in the presence of 
the \bc\ satisfy the Ward identities for \abar\ and accordingly can be 
expressed through chiral blocks for \abar. There are then two problems to be 
addressed. First, for a given CFT model whose chiral \alg\ contains \abar\ 
one wishes to construct all \bc s that preserve (at least) \abar, as well 
as arbitrary \corfu s on surfaces with such \bc s. Second, for any given 
\voa\ \abar\ one wants to know all CFT models that possess a chiral \alg\ 
\aext\ containing \abar. (In the context of \bc s, the task is actually 
quite the opposite, namely to determine all vertex operator sub\alg s of a 
given \voa\ $\aext\,{\equiv}\,\abulk$.) These two problems are
closely related, and they should be treated simultaneously. Partial answers
to both questions are known (see e.g.\ \cite{gann16} for the \class\ 
problem, \cite{fuSc112} for various aspects of \bc s, and \cite{fffs3} 
for the construction of \corfu s), but a general solution is still to be 
found. In general the number of conformally invariant \bc s 
is infinite. That one deals with only finitely many \bc s happens only 
in rather special circumstances where the vertex algebras involved are 
all rational. We emphasize that no such finiteness assumption enters our 
formalism; in particular we do not require semisimplicity of the relevant 
tensor categories.  Still, special situations with only finitely many 
boundary conditions can be of particular interest. An example is provided 
by the unitary models with $c\,{<}\,1$, for which
every irreducible \aext-module is finitely fully reducible as a
Virasoro module; in this case all conformal \bc s are known explicitly 
\cite{afos,fuSc9,bppz2}. More generally, in any rational CFT one may wish to 
restrict one's attention to the subclass of `rational' \bc s, for which by
definition the preserved sub\alg\ \abar\ is a rational \voa, such that every
irreducible \aext-module is finitely fully reducible as an \abar-module.  

In the present paper we investigate a proposal for a universal and model 
independent construction of all conformal \bc s. In fact, the construction 
encodes at the same time both a modular invariant for \abar\ \cite{boev123}, 
i.e.\ a `CFT in the bulk', and the conformal \bc s for that CFT. When combined 
with ideas from \cite{fffs2,fffs3}, it also opens the way to determine all 
\corfu s. And it sheds new light on the \class\ problem (both of modular 
invariants and of \bc s) as well; given a chiral algebra and hence an 
associated category \calc, what is to be classified are the `good' algebra 
objects in \calc.

We should motivate why in the last statement we formulate the task in
terms of a (representation) category rather than in purely algebraic terms.
Vertex algebras constitute, at least in our opinion, a quite intricate
mathematical structure. It can therefore be fruitful not to work directly
with such an algebra, but rather with its category \calcm\ of representations.
Indeed, much of the algebraic structure of a \voa\ \abar\ leads directly 
to interesting structure on its representation category. For instance, a 
chiral algebra is expected to possess coproduct-like structures, see 
\cite{fefk3} for an early discussion and \cite{huleX,huan3} for an
approach in the context of vertex algebras. As the representation category of
a chiral algebra, \calcm\ is therefore expected to carry the structure of
a tensor category. For rational \voa s \calcm\ is expected\,%
 \footnote{~The properties which equip \calcm\ with the structure of a \mtc\ 
 can be encoded in so-called Moore\hy Seiberg data \cite{mose3,fefk3,BAki}.
 It has actually not yet been proven that the \rep\ category of every rational
 VOA indeed possesses all features of a modular tensor category. But this
 property has been established for several classes of VOAs, compare e.g.\
 \cite{huan7,hule5}, and it is commonly expected that possible exceptions 
 should better be accounted for by an appropriate refinement of the
 qualification `rational'.} 
to be even modular. Another advantage of this approach is the following. 
Vertex algebras constitute just one particular mathematical formalization
of the physical idea of a chiral algebra. Other approaches exist which use 
different algebraic structures (see e.g.\ \cite{EVka,muge7}). They lead to 
categories that are, sometimes, equivalent to categories 
of representations of vertex algebras. A treatment on the level of
category theory therefore allows to deal with aspects of conformal field
theory that are common to all known formalizations of chiral algebras.
Yet another advantage is that similar categories also arise from other 
algebraic objects like quantum groups. Those objects are sometimes easier 
to understand than vertex algebras and can serve as a source of inspiration
for general constructions. (For instance, one can hope that non-compact
quantum groups will improve our understanding of certain non-compact
conformal field theories \cite{pote}.)

Once one adopts this point of view, one is lead to the following proposal that
is based on the observations in \cite{fuSc14,scfu3}, which in turn were
motivated by a comparison
of concrete information about \bc s in specific models with results 
from subfactor theory \cite{lore,loro,xu3,boev123,boek12} and about
the modularisation of premodular categories \cite{brug2,muge6}.  
The basic idea is to realize the elementary \bc s preserving \abar\ as the 
(absolutely) simple objects of a suitable category, which we like to call the 
{\em boundary category\/}. Non-simple objects of the boundary category 
correspond to `composite' \bc s. These are of particular interest in string 
theory, where the decomposition of a semisimple object into simple subobjects
specifies the so-called Chan\hy Paton multiplicities, which in turn are a 
source of gauge symmetries in space-time. According to \cite{rers} they also 
appear in the study of renormalization group flows on boundary conditions, 
even when the starting point corresponds to a simple boundary condition. 
(A category theoretic approach to certain specific boundary conditions 
in superconformal field theories has also been put forward in \cite{doug9}. 
It is in many respects different from our proposal; roughly speaking, our 
categories relate to the fusion ring of a CFT in the same manner as
the categories in \cite{doug9} should relate to the `chiral ring' of a 
superconformal field theory.) 

The boundary category is a tensor category with certain additional properties.  
It is constructed as a category of \A-modules in an underlying tensor category
\calc, with \A\ a suitable algebra object in \calc. An important technical 
point of our approach is that we try to minimize the assumptions about the 
category \calc. Indeed, upon a judicious choice of the relevant class of 
algebra objects, for some structure to be present in the resulting boundary 
category only that particular structure needs to be assumed for \calc.
For instance we do {\em not\/} have to require \calc\ to be modular, nor to
be semisimple. 
As is well known \cite{dolm3}, semisimplicity of the representation category  
of a vertex algebra amounts to rationality of the theory in the sense that 
there are only finitely many inequivalent simple modules. Thus our framework
offers some elements that will be useful to go beyond rational conformal 
field theories, and beyond boundary conditions that preserve a rational 
subalgebra in rational theories.  

In favorable situations the Grothendieck group of 
the boundary category inherits a ring structure from the tensor product. 
(This is certainly the case when the tensor bifunctor $\otimes$ is exact.)
When in addition the algebra object \A\ is commutative 
and has Frobenius\hy Schur indicator equal to one,
then the structure constants of the Grothendieck ring -- the 
{\sl fusion rules\/} -- coincide with the annulus coefficients of the CFT,
i.e.\ with the coefficients in an expansion of the \parfu\ on the annulus 
\wrtt character-like functions that are associated to the simple objects of 
the boundary category. Furthermore, the (quantum) dimension of an object of 
the boundary category gives us the `entropy' of the corresponding \bc, and
the $6j$-symbols of the boundary category provide the operator product 
coefficients of the so-called boundary fields. In the present contribution 
we will not, however, be concerned with this type of questions, but rather 
deal with category theoretic aspects of the proposal. Accordingly we also 
include various observations that are not immediately related to the 
construction of the boundary category but are of independent interest. A key 
role will be played by the concept of {\sl\hapl\ Frobenius \alg s\/} in 
abelian tensor categories, which is motivated by considerations in \cft;
the Frobenius property encodes s-t-duality of four-point chiral blocks in
the vacuum sector (compare figure \erf pp below), while \hapy\ corresponds 
to the uniqueness of the vacuum state. (These algebras should also possess 
a few further properties, which implement physical requirements, too. For 
details, see the main text.) A conjectural picture for the moduli space of 
\cfts\ of given Virasoro central charge is then that over each
point of the moduli space there lies a \hapl\ Frobenius algebra.
This picture is somewhat reminiscent of the role of Frobenius manifolds
in the study of integrable systems \cite{Dubr}.

Our category theoretic results are contained in sections 2 to 5.
Readers less interested in the mathematical development are invited to start
a first reading in section 6, where we summarize the salient features that 
are relevant in the context of conformal \bc s.
Afterwards we present, in sections 7 and 8, some details for specific 
situations of particular interest. We conclude by pointing out the 
relation of our construction with non-negative integral matrix \rep s 
(NIM-reps) of the fusion rules, obtaining in particular a criterion for
a NIM-rep to correspond to consistent \bc s 
for a torus partition function of extension type (section 9), 
and speculate about its interpretation in terms of \voa s (section 10).

\medskip\noindent
{\bf Note added}:\\
As mentioned below, ideas similar to some of ours have been expressed
a little earlier by A.A.\ Kirillov and V.\ Ostrik in \cite{kios,kirI14}
(compare also \cite{ostr} for a continuation of their work). 
After completing this contribution, other related work has come to our 
attention. Pioneering results on algebras in braided categories were 
obtained by B.\ Pareigis in \cite{pareX,pare23}.
A.\ Wassermann \cite{wassermann-unpublished} has considered,
in the context of quantum subgroups, (braided-)\,commutative algebras
in *-tensor categories. G.\ Moore \cite{moore-unpublished}
and G.\ Segal \cite{segal-unpublished,sega13} have developped an algebraic 
approach to \twodim\ {\em topological\/} field theories on surfaces with 
boundary that is closely related to the programme we have meanwhile pursued 
with I.\ Runkel \cite{fuRs,fuRs4} for conformal field theories; roughly
speaking, their results correspond to specializing the constructions in 
\cite{fuRs,fuRs4} to the modular tensor category of finite-dimensional 
complex vector spaces.

\section{Algebras in tensor categories}

In order to construct the boundary category which describes boundary 
conditions based on the \voa\ \abar, as well as for encoding the modular 
invariant torus partition function besides the \rep\ category $\calc\eq
\calcm$, one other basic datum is needed: a suitable algebra object \A\ 
in \calc. The boundary category is then the category of \A-modules in
\calc, and will accordingly be denoted by \calca. Such structures appear 
naturally in categories. In the context of vertex algebras, algebra objects 
are e.g.\ expected to arise in the following situation. Given a subalgebra 
\abar\ of a vertex algebra \aext, the ambient algebra \aext\ may be 
regarded as an \abar-module and thereby corresponds to an object $A$ in the 
tensor category \calcm. It can be expected (see e.g.\ \cite{kios}) that the
`product' $v\,{\star_z}w\,{:=}\,Y(v,z)w$ on the vertex algebra \aext\ induces
a product on the object \A\ that is commutative in the appropriate sense.

Eventually our desired construction will also involve module categories in
categories \calc\ that are closely related to, but different from, \calcm.
Therefore in the sequel we do not restrict our attention to categories that
are modular. Rather we consider, more generally, strict tensor\,%
 \footnote{~As a duality in a category is a rather subtle structure, we do
 not include its existence in the definition of the term `tensor category'.
 Thus what we refer to as a tensor category is often called a {\sl monoidal\/}
 category in the literature. Dualities will be studied in section 3.}
categories which are {\sl abelian\/} and for which the tensor product 
of morphisms is $k$-{\sl additive\/} in each factor, where $k$ is the 
ground ring. For instance, \calc\ may be a boundary category itself; 
boundary categories do not possess, in general, a braiding in the usual 
sense, and their Grothendieck ring can be non-commutative. We also
do not require the category \calc\ to be semisimple. Such a requirement
would inevitably restrict us to categories that correspond to rational \cfts, 
since rationality of a vertex algebra amounts to requiring semisimplicity 
of its representation category \cite{dolm3}.  As we do not assume 
semisimplicity, our setting is more general than the one in 
\cite{kios,kirI14}. Nevertheless, some of our arguments (like the use of 
induction and restriction functors) and results (in particular those which 
follow when \calc\ {\em is\/} semisimple, such as the decomposition 
\erf ab of $A\Oti A$ and
proposition \ref{25}) are parallel to the reasoning in \cite{kios,kirI14}.

We use the following basic notations and concepts. The class of objects 
of \calc\ is denoted by $\obj(\calc)$, and the morphisms from
$X\iN\obj(\calc)$ to $Y\iN\obj(\calc)$ by $\Hom(X,Y)$. The endomorphisms 
$\Hom(\one,\one)\,{=:}\,k$ of the tensor unit $\one\iN\obj(\calc)$ form a 
commutative associative unital ring $k$, called the ground ring. We refer 
to the elements of $k$ as {\sl scalars\/} and denote the unit $\id_\smallone$ 
of $k$ by 1. Each morphism set $\Hom(X,Y)$ is a bimodule over $k$, with 
commuting left and right actions of $k$.

Throughout the rest of the paper \calc\ will be assumed to be a strict abelian 
\tc\ for which the tensor product is $k$-additive in each factor.  
As usual, owing to the coherence theorems the assumption of strictness does 
not restrict the validity of our results; it enables us to suppress
the associativity and unit constraints, the inclusion of which would make 
many formulas more difficult to read.

                    \medskip
                    \vfill

We start with the following fundamental definitions.

\definition
(i)\ An {\sl algebra\/} in $\calc$ is a triple $(A,m,\eta)$, where 
$A\iN\obj(\calc)$, $m\iN\Hom(A\Oti A,A)$ and $\eta\iN\Hom(\one,A)$,
obeying
  \be  m \circ (\id_A\oti m) = m\circ (m\oti\id_A) \labl as \end{equation} 
and 
  \be  m \circ (\id_A\oti \eta) = \id_A =
  m \circ(\eta\oti \id_A) \,.  \labl un  \end{equation}

                    \vfill
\vskip3pt\noindent
(ii)\ A {\sl co-algebra} in $\calc$ is a triple $(C,\Delta,\eps)$, where 
$C\iN\obj(\calc)$, $\Delta\iN\Hom(C,C\Oti C)$ and $\eps\iN\Hom(C,\one)$, obeying
  \be (\Delta\oti\id_C)\circ\Delta
  = (\id_C\oti\Delta)\circ\Delta   \end{equation}
and
  \be (\id_C\oti\eps) \circ \Delta = \id_C
  = (\eps\oti\id_C) \circ\Delta \,.  \end{equation}

                    \vfill
\vskip3pt\noindent
(iii)\ Given an isomorphism $\tilde c\iN\Hom(A\Oti A,A\Oti A)$ satisfying
  \be \bearl  (m\oti\id_A) \circ (\id_A\oti\tilde c) \circ
  (\tilde c\oti\id_A) = \tilde c \circ (\id_A\oti m) \,, \\{}\\[-.6em]
  (\id_A\oti m) \circ (\tilde c\oti\id_A) \circ
  (\id_A\oti\tilde c)  = \tilde c \circ (m\oti\id_A)  \eear
  \end{equation}
                      \newpage
and
  \be  \tilde c \circ (\id_A\oti\eta) = \eta \oti \id_A \,,  \qquad
  \tilde c \circ (\eta\oti\id_A) = \id_A\oti\eta \,,  \end{equation}
an algebra \A\ in \calc\ is called $\tilde c$-{\sl commutative\/} iff
  \be  m \circ \tilde c = m \,, \end{equation}
and a co-algebra $C$ in \calc\ is called $\tilde c$-{\sl cocommutative\/} iff
  \be  \tilde c \circ \Delta = \Delta \,.  \end{equation}

\vskip.3em
                    \vfill

\remark
1.\ The morphism $m$ is called multiplication or product of \A, and $\eta$ is 
called the unit morphism of \A, or unit, for short. 
$\Delta$ and $\eps$ are called the co\-pro\-duct and counit of $C$, \resp,
and the property \erf as is referred to as the associativity axiom, etc.
\\
2.\ The isomorphism $\tilde c$ used in the definition of
$\tilde c$-(co-)commutativity is a baby version of a {\sl\swap\/}, see 
definition \ref{sw} below. The notion of a \swap\ is due to 
A.\ Bru\-gu\-i\`e\-res \cite{brug+}.
\\
3.\ When \calc\ is the category of \findim\ vector spaces over a field $F$, 
a \mbox{(co-)al}\-ge\-bra in \calc\ is an $F$-(co-)\alg\ in the ordinary sense.
\\
4.\ 
When \A\ is both an \alg\ and a co-\alg, then the fact that the endomorphisms
of the tensor unit $\one$ are just the ground ring implies 
in particular that $\eps\cir\eta\eq\gam\,\id_\one$ for some $\gam\iN k$.
Also note that $f\,{:=}\,\eta\cir\eps$ is non-zero (because of
$m\cir(f\Oti\id_A)\cir\Delta\eq\id_A$, which follows by the unit and counit 
properties) and satisfies $f{\circ}f\eq\gam f$. 
\\
5.\ Besides with the concatenation, the endomorphisms $\Hom(A,A)$ of
an object \A\ that is both an \alg\ and a co-\alg\
can also be endowed with another product structure, via
  \be  f * g := m \circ (f\oti g) \circ \Delta \,.  \labl fg \end{equation}
Due to (co-)associativity this operation is associative, and due to the
(co-)unit axioms it has the endomorphism $\eta\cir\eps$ as a unit.

\vskip.3em
                    \vfill

In the sequel we reserve the symbol $A$ for an object of \calc\ that is 
both an \alg\ and a co-\alg, but for brevity still just use the term algebra.  
For the applications we have in mind we need to further restrict our attention 
to \alg s \A\ that possess some additional structure. More concretely,
it makes sense to impose the following compatibility
between the algebra and co-algebra structure:

                    \vfill
\definition
(i)\ A {\sl Frobenius algebra\/} in a \tc\ $\calc$ is a quintuple
$(A,m,\eta,\Delta,\eps)$ such that $(A,m,\eta)$ is an algebra in \calc, 
$(A,\Delta,\eps)$ a co-algebra in \calc, and
  \be  (\id_A\oti m) \circ (\Delta\oti\id_A)
  = \Delta \circ m = (m\oti\id_A) \circ (\id_A\oti\Delta) \,.  \labl1f
  \end{equation}
(ii)\ A {\sl special algebra\/} in an abelian \tc\ \calc\ is a quintuple
$(A,m,\eta,\Delta,\eps)$ such that $(A,m,\eta)$ is an algebra in \calc,
 $(A,\Delta,\eps)$ is a co-algebra in \calc, and 
  \be  m \circ \Delta = \bet\, \id_A  \labl3f  \end{equation}
as well as
  \be  \eps \circ \eta = \gam\, \id_\smallone \equiv \gam
  \labl2f  \end{equation}
for suitable invertible scalars $\bet,\gam\iN k$.

          \newpage 

\remark \labl st
1.\ As already mentioned, the equality \erf2f is in fact fulfilled 
automatically; but we do have to require that $\gam$ be {\em invertible\/}.
Also note that it already follows from the definition of the ground ring that
  \be  \eps \circ m \circ \Delta \circ \eta = \beta  \labl1b  \end{equation}
with $\beta\iN k$.
\\
2.\ Property \erf3f just says that $\bet^{-1}\Delta$ is a right-inverse of 
$m$. This, in turn, tells us that $A$ is a subobject of $A\Oti A$. Similarly,
property \erf2f implies that $\one$ is a distinguished subobject of \A; 
later on we  will consider a class of algebra objects for which in addition 
$\Hom(\one,A)\,{\cong}\,k$, see formula \erf1i.
\\
3.\ In any \tc, the tensor unit $\one$ is a special Frobenius \alg\ with
$\bet\,{\equiv}\,\gam\eq1$. The product is just given by the unit constraint
$\ell_\smallone\,{\equiv}\,r_\smallone$, the coproduct is its inverse,
and the unit and counit are given by $\id_\smallone$.
\\
4.\ For any finite set M, the set ${\mathcal F}{\rm M}$ of complex functions
on M is a commutative \alg\ in the category of \findim\ complex vector
spaces, with the product given by point-wise multiplication and the unit being
the constant function 1. In terms of a basis $\{\delta_{\rm g}\,|\,{\rm g}
\iN{\rm M}\}$ of ${\mathcal F}{\rm M}$ (with $\delta_{\rm g}({\rm h})\eq
\delta_{\rm g,h}$ for ${\rm g,h}\iN{\rm M}$), the product reads
$m(\delta_{\rm g},\delta_{\rm h})\eq\delta_{\rm g,h}\,\delta_{\rm g}$,
and the unit is $1\eq\sum_{{\rm g}\in{\rm M}}\delta_{\rm g}$. By setting
  \be  \Delta(\delta_{\rm g}) := \delta_{\rm g}\oti\delta_{\rm g}
  \qquad{\rm and}\qquad  \eps(\delta_{\rm g}) := 1  \end{equation}
for all ${\rm g}\iN{\rm M}$ -- that is, the coproduct is the diagonal and
the counit the integral -- ${\mathcal F}{\rm M}$ becomes a special
Frobenius \alg\ with $\bet\eq1$ and $\gam\eq|{\rm M}|$.\,%
 \footnote{~If M is a group, then the algebra ${\mathcal F}{\rm M}$ admits
 another co-\alg\ structure, which can be extended to a Hopf \alg\ rather
 than a Frobenius algebra, on the same vector space. This latter coproduct
 implements the group multiplication, i.e.\ $\tilde\Delta(\delta_{\rm g})
 \eq\sum_{{\rm g}\in{\rm M}}\delta_{\rm h}\Oti\delta_{\rm gh^{-1}}$.}
A Frobenius \alg\ equivalent to ${\mathcal F}{\rm G}$, where G is
a finite group, is present \cite{kirI14} in the modular tensor category
for a so-called holomorphic orbifold \cft\ with orbifold group G.
\\
5.\ For symmetric tensor categories a definition of Frobenius object similar 
to the one used here has been given in \cite{stri2}.\,%
   \footnote{~We are grateful to M.\ M\"uger for drawing our attention to
   the concept of Frobenius algebras in tensor categories, and to the
   references \cite{stri2,loro}.}
The motivation in \cite{stri2} is that the 
spectrum of any finite groupoid provides an example of such an object.
In the category of vector spaces, one usually defines a Frobenius \alg\ 
as an \alg\ with a non-degenerate invariant pairing; in the present setting 
an analogue of such a pairing is supplied by the morphism $\eps\cir m$.
\\
6.\ In the framework of *-categories, a special Frobenius algebra is known 
as a {\sl Q-sys\-tem\/} \cite{long6,loro}. In that case the product and 
coproduct, and the unit and counit, respectively, are *'s of each other, and 
the Frobenius axiom \erf1f can be derived from the other axioms.
\\
7.\ In the vertex operator setting, where the algebra object corresponds
to a pair consisting of a vertex algebra \aext\ and a subalgebra $\abar\,
{\subset}\,\aext$, the Frobenius property is related to the fact that 
four-point conformal blocks for the vacuum on the sphere can be factorized 
in two different ways. In physics terminology, the factorizations are called 
the s-channel and t-channel, \resp, and the relationship between the two is 
known as s-t-duality. Pictorially,
  \be \begin{picture}(200,52)(68,38) 
  \put(17,7){\begin{picture}(0,0)(0,0)
  \scalebox{.16}{\includegraphics{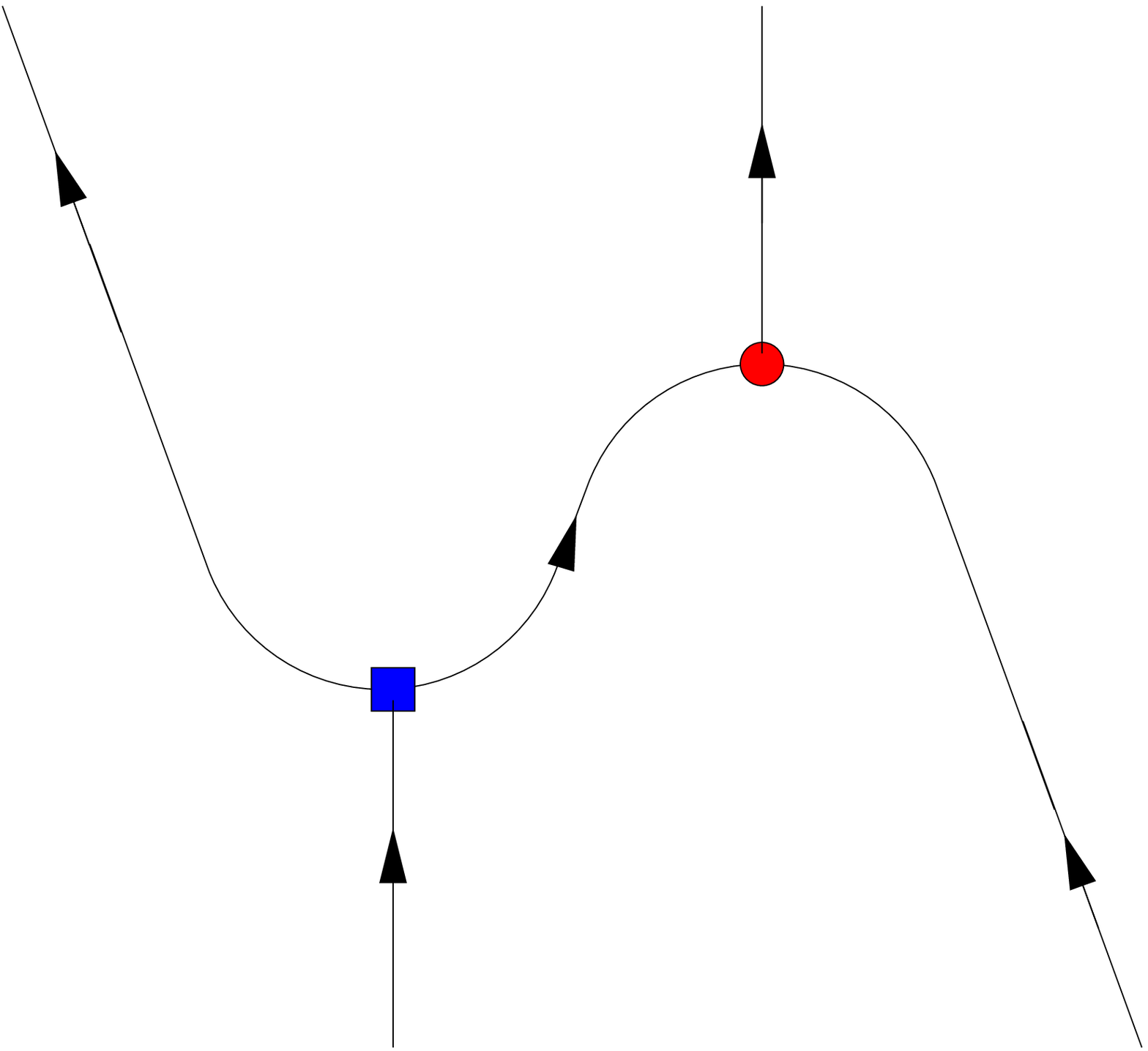}} \end{picture}}
  \put(109,39){$=$}
  \put(133,0){\begin{picture}(0,0)(0,0)
  \scalebox{.16}{\includegraphics{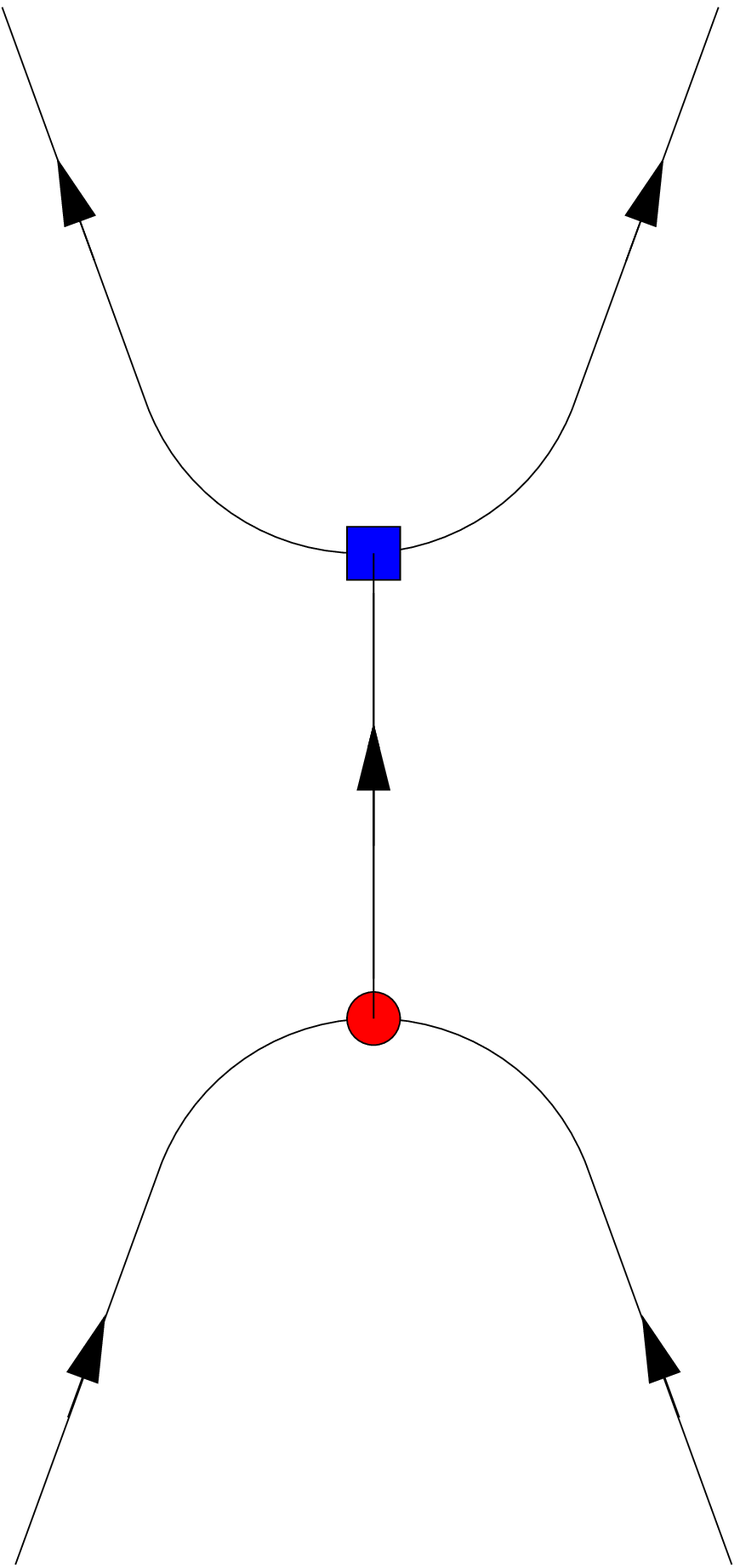}} \end{picture}}
  \put(191,39){$=$}
  \put(215,7){\begin{picture}(0,0)(0,0)
  \scalebox{.16}{\includegraphics{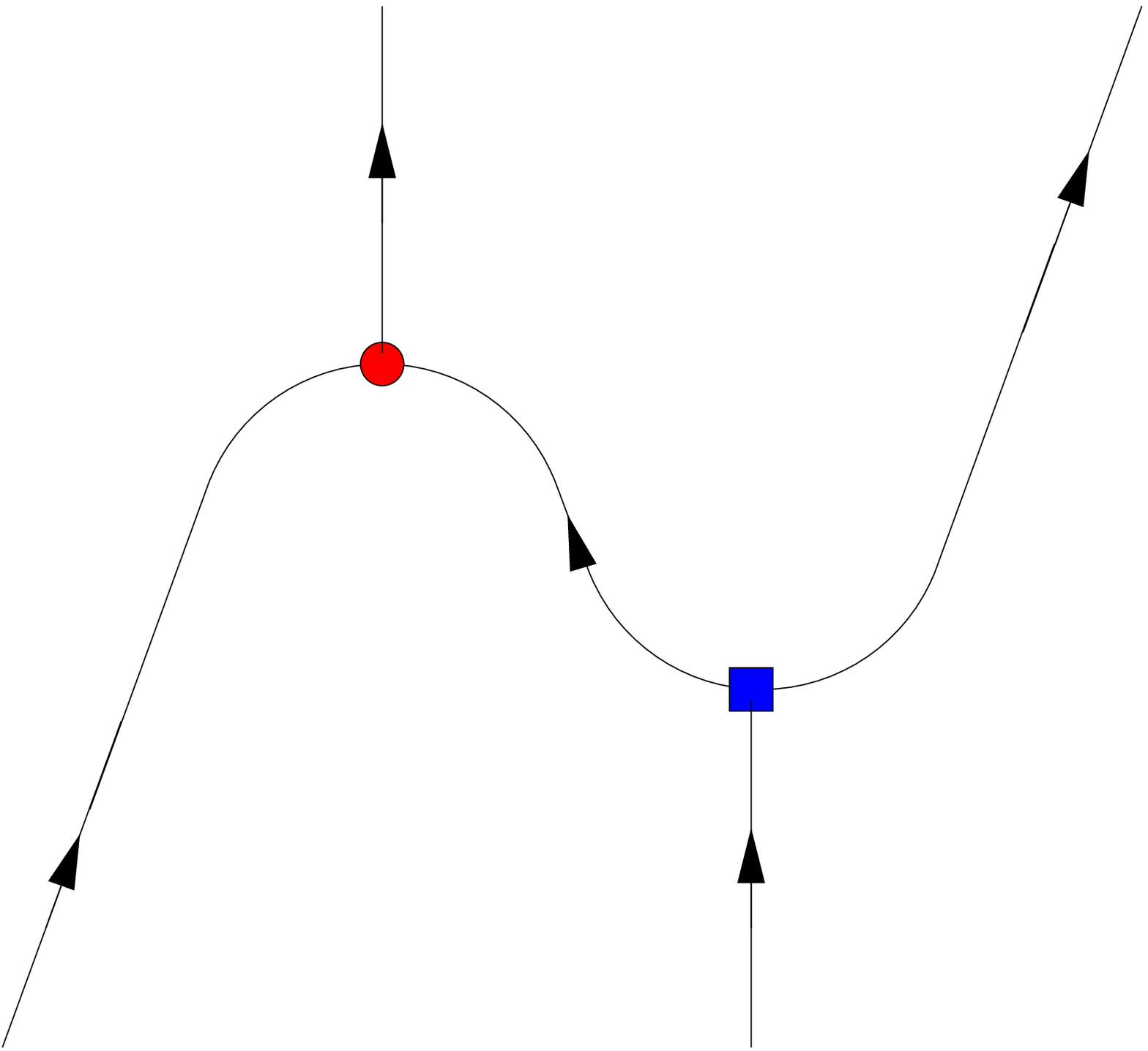}} \end{picture}}
  \end{picture} \labl pp  \end{equation}

\vskip3.9em 
\noindent
In this respect, the Frobenius property is thus quite 
similar to associativity, which may be regarded as encoding s-t-duality for 
a situation with three `incoming' and one `outgoing' state.

\vskip.53em

\lemma
For a \tcc\ Frobenius algebra one has
  \be   (m\oti m) \circ (\id_A\oti\tilde c\oti\id_A)
  \circ (\Delta\oti\Delta) = \Delta \circ m \circ \Delta \circ m \,.
  \end{equation}
In particular, when \A\ is special Frobenius, then
  \be   (m\oti m) \circ (\id_A\oti\tilde c\oti\id_A)
  \circ (\Delta\oti\Delta) = \bet\, \Delta \circ m \,.  \labl bi
  \end{equation}
\proof\dd
By elementary use of the axioms. As an illustration of the type of
manipulations that are involved, we nevertheless present 
the calculation explicitly:
  \be 
  \hsp{-.3}
  \bearll  (m\oti m) \circ (\id_A\oti\tilde c\oti\id_A)
  \circ (\Delta\oti\Delta)
  \\{}\\[-.5em] \hsp{1.3}
  = (m\oti m) \circ (\id_A\oti\Delta\oti\id_A)
  \circ (\tilde c\oti\id_A) \circ (\id_A\oti\Delta)
  & \mbox{\footnotesize by $\tilde c$-commutativity}
    \\[-.2em]&\mbox{\footnotesize and functoriality of $\tilde c$}
  \\{}\\[-.88em] \hsp{1.3}
  = (m\oti\id_A) \cir (\id_A\oti\Delta) \cir (\id_A\oti m)
  \cir (\tilde c\oti\id_A) \cir (\id_A\oti\Delta)\!\!\!\!%
  & \mbox{\footnotesize by Frobenius property}
  \\{}\\[-.6em]
   \hsp{1.3}
  = \Delta \circ m \circ (\id_A\oti m)
  \circ (\tilde c\oti\id_A) \circ (\id_A\oti\Delta)
  & \mbox{\footnotesize by Frobenius property}
  \\{}\\[-.6em] \hsp{1.3}
  = \Delta \circ m \circ (m\oti\id_A)
  \circ (\tilde c\oti\id_A) \circ (\id_A\oti\Delta)
  & \mbox{\footnotesize by associativity}
  \\{}\\[-.6em] \hsp{1.3}
  = \Delta \circ m \circ (m\oti\id_A)\circ (\id_A\oti\Delta)
  & \mbox{\footnotesize by $\tilde c$-commutativity}
  \\{}\\[-.6em] \hsp{1.3}
  = \Delta \circ m \circ \Delta \circ m
  & \mbox{\footnotesize by Frobenius property.}
\\[-.75em] {}
  \eear \end{equation}
{}\\[-.9em]\phantom.\hfill\checkmark

\medskip

\noindent
Note that formula \erf bi is a projective version of the defining property
of a bi-\alg\ in \calc.

\medskip

\section{Dualities and the FS indicator}

\definition \labl du
(i)\ A (right) {\sl duality\/} on a tensor category \calc\ associates to 
every $X\iN\obj(\calc)$ an object $X\Vee\iN\obj(\calc)$ and morphisms
  \be  b_X \in\Hom(\one,X\Oti X\Vee)\,, \qquad
  d_X \in \Hom(X\Vee\Oti X,\one)  \labl1d  \end{equation}
such that
  \be  (\id_X\oti d_X) \cir (b_X\oti\id_X) = \id_X \,,
  \quad\  (d_X\oti \id_{X\Vee}) \cir (\id_{X\Vee}\oti b_X)
  = \id_{X\Vee}  \,,  \labl2d  \end{equation}
and to every morphism $f\iN\Hom(X,Y)$ the morphism
  \be  f\Vee := (d_Y\oti\id_{X\Vee}) \circ (\id_{Y\Vee}\oti f\oti\id_{X\Vee})
  \circ (\id_{Y\Vee}\oti b_X) \ \in \Hom(Y\Vee\!{,}\,X\Vee) \,.
  \end{equation}
(ii)\ $X\Vee$ is called the object (right-) {\sl dual\/} to the object $X$,
and $f\Vee$ the morphism (right-) {\sl dual\/} to the morphism $f$.
\\[.32em]
(iii)\ When $X\Vee$ is isomorphic to $X$, then $X$ is called (right-) 
{\sl self-dual\/}.

\medskip

\remark
1.\ Up to equivalence of categories, the duality is unique. In more detail, 
the dual object $X\Vee$ of $X$ is unique up to a canonical isomorphism, and 
using those isomorphisms one shows that the duality morphisms are unique up
to canonical isomorphism as well; these isomorphisms furnish a natural 
transformation.
\\[.19em]
2.\ The tensor unit is self-dual, $\one\Vee{\cong}\,\one$ (isomorphisms are
obtained by composing duality morphisms with unit constraints);
but we will not assume that $\one\Vee{=}\,\one$.
Also, the double dual $(X\Vee)\Vee_{}$ 
need not coincide with, nor even be isomorphic to, $X$. 
\\[.19em]
3.\ It is not required that $b$ and $d$ satisfy any compatibility 
conditions with other structures.  Thus for instance in modular categories 
the compatibility $(\theta_X\Oti\id_{X\Vee})$\linebreak[0]$\cir b_X
$\linebreak[0]$\eq(\id_X\Oti\,\theta_{X\Vee})\cir b_X$ of the duality 
with the twist endomorphisms must be imposed as a separate axiom.
On the other hand, it follows immediately that
  \be  (f\oti\id_{X\Vee}) \circ b_X = (\id_Y \otimes f\Vee) \circ b_Y
  \,.  \end{equation}
for all $X,Y\iN\obj(\calc)$ and all $f\iN\Hom(X,Y)$. 
(Thus e.g.\ the compatibility axiom of a modular category just mentioned 
is nothing else than $(\theta_X)\Vee{=}\,\theta_{X\Vee}$.)
Also, in full generality one has $(f\cir g)\Vee\eq
g\Vee{\circ}\, f\Vee$ for any two composable morphisms $f,g$, as well as
$(\id_X)\Vee\eq\id_{X\Vee}$ for every $X\iN\obj(\calc)$.
\\[.19em]
4.\ Further, $(X\Oti Y)\Vee$ is isomorphic to $Y\Vee\Oti X\Vee$, and
for any three objects $X,Y,Z$ there are natural bijections
  \be  \hsp{-.42} \Hom(X\Oti Y,Z) \cong \Hom(X,Z\Oti Y\Vee) \;\ {\rm and}\;\ 
  \Hom(X,Y{\otimes}\,Z) \cong \Hom(Y\Vee\!{\otimes}X,Z)  \labl34
  \end{equation}
of $k$-bimodules.
These properties imply that the operation of replacing $X$ by $X\Vee$
and $f$ by $f\Vee$ furnishes a functor from \calc\ to $\calc^{\rm op}$ and 
moreover, that the functor of tensoring from the right with $X$ is left-adjoint 
to tensoring from the right with $X\Vee$, while the functor of tensoring from 
the left with $X\Vee$ is left-adjoint to tensoring from the left with $X$
(see e.g.\ \cite[XIV.2.2]{KAss}). It follows that in a category with right 
duality the functor defined by tensoring with an object from the left is exact.
\\[.19em]
5.\ In a tensor category with duality, for every \alg\ \A\ there is a natural
co-\alg\ structure on $\Av$ (and vice versa), which is obtained by setting
  \be \Delta := m\Vee  \qquad{\rm and}\qquad  \eps := \eta\Vee \,.  \end{equation} 

\medskip

\lemma  \labl0f
A Frobenius \alg\ \A\ in a tensor category with duality is self-dual,
$A\,{\cong}\,\Av$. An isomorphism $\Phi\iN\Hom(A,\Av)$\, is given by
  \be  \Phi := (\eps\oti\id_\Av) \circ (m\oti\id_\Av)
  \circ (\id_A\oti b_A) \,.  \labl4f  \end{equation}
Moreover, if \A\ is \tcc, then it is \tccc\ as well.
\proof\dd
The morphism 
  \be  \Phi^- := (d_A\oti\id_A) \circ (\id_\Av\oti\Delta)
  \circ (\id_\Av\oti\eta) \ \in \Hom(\Av\!{,}\,A)  \labl5f
  \end{equation}
is inverse to $\Phi$. That it is a right-inverse follows by first
using the duality axiom \erf2d, then the Frobenius property \erf1f
and then the defining properties of the unit and counit. That it is
a left-inverse is seen by applying these identities in the reverse order.
\\
To prove the statement about $\tilde c$-cocom\-mu\-ta\-ti\-vity,
one composes $\tilde c\cir\Delta$ with the identity morphism written as 
$\id_A\Oti(\Phi^{-1}{\circ}\,\Phi)$
and then uses successively functoriality of $\tilde c$, \tcy, the
Frobenius property, the counit property, \tcy, the duality axiom,
once again the Frobenius property, and finally the unit property.
\hfill\checkmark

\medskip

\remark
An \alg\ \A\ in a braided tensor category with duality which is 
commutative (with respect to the braiding $c_{A,A}$), has
non-vanishing dimension, and for which \erf4f is an isomorphism, is
called a {\sl rigid\/} \calc-\alg\ in \cite{kios}.

\medskip

\definition
(i)\ An {\sl autonomous\/} 
category is a \tc\ \calc\ with a right duality \erf1d as well
as a {\sl left duality\/}. The latter consists of an association of a
$\Eev X\iN\obj(\calc)$ to every object $X$ of \calc, of morphisms
  \be \tilde b_X \in\Hom(\one,\Eev X\Oti X) \,, \qquad
  \tilde d_X \in \Hom(X\Oti\Eev X,\one)  \end{equation}
for every $X\iN\obj(\calc)$, satisfying
  \be  (\tilde d_X\oti\id_X)\cir(\id_X\oti\tilde b_X) = \id_X \,, \quad\
  (\id_{X\Vee}\oti\tilde d_X)\cir(\tilde b_X\oti\id_{X\Vee})
  = \id_{X\Vee} \,,  \labl3d  \end{equation}
and of
  \be  \Eev f := (\id_{\Eev X}\oti\tilde d_Y) \circ (\id_{\eev X}\oti f
  \oti\id_{\eev Y}) \circ (\tilde b_X\oti\id_{\eev Y})
  \ \in \Hom(\eev Y{,}\Eev X) \end{equation}
to every morphism $f\iN\Hom(X,Y)$.  
\\[.2em]
(ii)\ A {\sl sovereign\/} 
category is an autonomous category for which $\Eev X\eq X\Vee$ for all
$X\iN\obj(\calc)$ and $\Eev f\eq f\Vee$ for all $f\iN\Hom(X,Y)$ and all
$X,Y\iN\obj(\calc)$.
\\[.3em]
(iii)\ In a sovereign category, the left and right {\sl traces\/} of an
endomorphism $f\iN\Hom$\linebreak[0]$(X,X)$ are the scalars
  \be
  \trl(f) :=  d_X \circ (\id_{X\Vee}\oti f) \circ \tilde b_X \,, \qquad
  \trr(f) := \tilde d_X \circ (f\oti\id_{X\Vee}) \circ b_X \,,
  \end{equation}
and the left and right {\sl dimensions\/} of an object $X$ are the scalars
  \be  \diml(X) := \trl(\id_X) \,, \qquad \dimr(X) := \trr(\id_X) \,.
  \end{equation}
(iv) A {\sl spherical\/} category is a sovereign
category for which $\trl(f)\eq\trr(f)$ for all $f\iN\Hom(X,X)$ and all
$X\iN\obj(\calc)$.

\medskip

\remark
1.\ Analogously to \erf34, in the presence of a left duality we have 
isomorphisms
  \be  \hsp{-.3} \Hom(X\Oti Y,Z) \cong \Hom(Y,\Eev X\Oti Z) \,, \quad \
  \Hom(X,Y{\otimes}\,Z) \cong \Hom(X\Oti\Eev Z,Y) \,.  \end{equation}
Note that for a general autonomous category analogous statements with
$\Eev X$ replaced by $X\Vee$ etc.\ are not valid.
\\[.19em]
2.\ In an autonomous category one has $(\Eev X)\Vee\,{\cong}\,X\,{\cong}\,
\eev{(X\Vee)}$; hence in a sovereign category $(X\Vee)\Vee\,{\cong}\,X$
-- an isomorphism is e.g.\ provided by
$(\id_{(X\Vee)\Vee}\Oti d_X)\cir(\tilde b_{X\Vee}\Oti\id_X)$, with inverse
$(\id_X\Oti\tilde d_{X\Vee})\cir(b_X\Oti\id_{(X\Vee)\Vee})$ --
but still $(X\Vee)\Vee$ can be different from $X$. A notion similar 
to sovereignty is that of a {\sl pivotal\/} category
\cite{frye,bawe5}, where the primary structure is the isomorphism
between $X$ and $(X\Vee)\Vee$ rather than the left duality.
\\[.19em]
3.\ The traces in a sovereign category are cyclic, so that the term `trace'
is appropriate. This property depends crucially on the equality between left 
and right dual morphisms. Indeed, for every pair $f\iN\Hom(X,Y)$ and 
$g\iN\Hom(Y,X)$, by inserting the left duality axiom 
    $\id_X\eq(\tilde d_X\Oti\id_X)\cir(\id_X\Oti\tilde b_X)$,
\resp\ a right duality axiom for $\id_Y$, one shows that
  \be  \trl(f\cir g) = d_Y \circ (\Eev g\oti f) \circ \tilde b_X \quad\
  {\rm and}\quad\  \trl(g\cir f) = d_Y \circ (g\Vee\Oti\,f) \circ \tilde b_X
  \,.  \end{equation} 
Thus $\Eev g\eq g\Vee$ implies the cyclicity $\trl(f\cir g)\eq\trl(g\cir f)$;
analogously one shows that $\trr(f\cir g)\eq\trr(g\cir f)$.
\\
As a consequence of the cyclicity, traces depend only on isomorphism classes,
in the sense that for every isomorphism $g\iN\Hom(Y,X)$ and every endomorphism 
$f\iN\Hom(X,X)$ one has $\trl(gfg^{-1})\eq\trl(f)$ and $\trr(gfg^{-1})\eq
\trr(f)$. In particular, dimensions are constant on isomorphism classes of 
objects. Moreover, using the isomorphism between $(X\Vee)\Vee$ and $X$ 
one shows that $\trl(f)\eq\trr(f\Vee)$ and $\trr(f)\eq\trl(\Eev f)\,{\equiv}
\,\trl(f\Vee)$, and hence $\dimr(X)\eq\diml(X\Vee)$. It follows that
self-dual objects have equal left and right dimensions,
  \be  \diml(X) = \dimr(X) \quad\ {\rm if}\;\ X\,{\cong}\,X\Vee  \labl dd
  \,.  \end{equation}
In the case of self-dual objects we will therefore simply write $\dim(X)$; 
for instance, the (left and right) dimension of a Frobenius \alg\ \A\ will 
be denoted as $\dim(A)$.
\\[.19em]
4.\ The left (and analogously, the right) trace satisfies $\trl(f{+}g)
\eq\trl(f)\,{+}\,\trl(g)$ and $\trl(f\Oti g)\eq\trl(f)\,{\cdot}\,\trl(g)$.
In particular, left and right dimensions are additive under direct sums and
multiplicative under tensor products. 

               \vfill\noindent
5.\ In a category with left duality, tensoring with an object from the right 
yields an exact functor. Together with the corresponding statement for
a right duality, it follows that in an autonomous category tensoring with an 
object is two-sided exact.

               \vfill\noindent
6.\ Spherical categories allow the construction of invariants of links in
the three-sphe\-re (see e.g.\ \cite{bawe2}), which explains their name. The 
condition on the equality of traces implies that the invariant of a link in 
three dimensions does not depend on the two-dimensional projection of the link.
Now the calculation of \corfu s in CFT presumes the existence of invariants 
of certain three-manifolds. Still, we do not have to require that boundary 
categories are modular. The absence of this property should be compensated by 
additional structures on the three-manifold, which are similar to those of 
spin- \cite{blma,beli,sawi3} and so-called $\pi$- \cite{tura10} manifolds.

                \vfill
\medskip

\lemma \labl3r
Every object $A$ of a 
sovereign category \calc\ that can be written in the form
  \be  A = A_{X,X} := X \oti X\Vee  \labl ax  \end{equation}
is a special Frobenius algebra in $\calc$, with parameters
  \be  \bet = \diml(X)\,, \quad\ \gam = \dimr(X) \,.  \end{equation}

                \vfill

\proof\dd 
Define the product, coproduct, unit and counit by
  \be  \bearll
  m := \id_X \oti d_X \oti \id_{X\Vee} \,,\ \ & \eta := b_X \,, \\{}\\[-.5em]
  \Delta := \id_X \oti\tilde b_X\oti\id_{X\Vee} \,, & \eps := \tilde d_X \,.
  \eear \end{equation}
The statement then follows by straight-forward manipulations involving 
the defining properties of the duality morphisms.
\hfill\checkmark

                 \newpage 

\remark
1.\ This observation generalizes as follows. For any pair of objects $X,Y\iN
\obj(\calc)$ one can introduce the object $A_{X,Y}\,{:=}\,X\Oti Y\Vee$. There 
are then natural multiplications from $A_{X,Y}\Oti A_{Y,Z}$ to $A_{X,Z}$. 
Thereby every sovereign tensor category gives rise to an algebroid, in the 
category theoretic sense. In the case of the boundary category, it is this 
algebroid that furnishes the algebraic structure of the `open string sector'.

               \vfill\noindent
2.\ In the application of \cft\ to strings, the \alg\ structure on $A_{X,X}$   
also gives rise to gauge groups in open string theory. The type I
string in 10 flat dimensions, for example, possesses a single Neumann boundary 
condition, known as the D9-brane. The corresponding boundary category
is therefore equivalent to the category of vector spaces. 
Having a `Chan-Paton multiplicity' $n$ for the D9-brane corresponds to
considering an object $X$ that is the direct sum of $n$ simple objects, i.e.\
an $n$-dimensional vector space. $A_{X,X}$ can then be identified with the 
matrix algebra of $n{\times}n$\,-matrices. The so-called orientifold projection
reduces this to an algebra of antisymmetric matrices. With the value $n\eq32$
that results from the tadpole cancellation condition, one thus arrives at the 
Lie algebra of the usual gauge group Spin$(32)/\zet_2$.

               \vfill\noindent
3.\ Objects of the form \erf ax are in particular `canonical objects'
in the sense of \cite{loro}.

\vskip.75em
               \vfill

For any Frobenius algebra in a sovereign category \calc\ there is
another distinguished isomorphism in $\Hom(A,\Av)$ besides \erf4f, namely
  \be \bearll \tilde\Phi \!\!\!
  &{:}{=}\; (\id_\Av\oti d_A) \circ (\id_\Av\oti\Phi\oti\id_A) \circ
  (\tilde b_A \oti\id_A)
  \\{}\\[-.7em] 
  &\;= (\id_\Av\oti\eps) \circ (\id_\Av\oti m) \circ (\tilde b_A\oti\id_A)
  \,, \eear  \labl ft  \end{equation}
with inverse $\tilde\Phi^{-1}\eq(\id_A\oti\tilde d_A)\circ(\Delta\oti\id_\Av)
\circ(\eta\oti\id_\Av)$.
               \vfill

We can also `invert' the formulas for $\Phi$ and $\tilde\Phi$ and their
inverses so as to obtain the following expressions for the duality morphisms
of $A$:
  \be \bearll 
  b_A = (\id_A \oti \Phi) \circ \Delta \circ \eta \,,  & 
  \tilde b_A = (\tilde \Phi \oti \id_A) \circ \Delta \circ \eta \,, 
  \\{}\\[-.5em]
  d_A = \eps \circ m \circ (\Phi^{-1}_{} \oti \id_A) \,,\ &
  \tilde d_A = \eps \circ m \circ (\id_A \oti \tilde\Phi^{-1}_{}) \,.
  \eear  \end{equation} 
The same kind of calculations show that $m\cir(\id_A\Oti\Phi^{-1}_{})
\cir b_A$, which must be a multiple of $\eta$, is in fact
  \be  m \circ (\id_A \oti \Phi^{-1}_{}) \circ b_A = \beta\, \eta \,,
  \end{equation}
where $\beta$ is the scalar introduced in formula \erf1b, and analogously
  \be \bearl   m \circ (\tilde\Phi^{-1}_{} \oti \id_A) \circ \tilde b_A
  = \beta\, \eta \,, \\{}\\[-.7em]
  d_A\circ (\Phi \oti \id_A) \circ \Delta = \beta\, \eps
  = \tilde d_A\circ (\id_A \oti \tilde\Phi) \circ \Delta \,. 
  \eear  \end{equation}
Further, when \erf2f holds with invertible $\gam$, then we also have
  \be \bearl
  (d_A\oti\id_A)\cir(\id_\Av\oti\Delta)\cir\tilde b_A = (\dimL(A)/\gam)\,\eta\,,
  \\{}\\[-.5em]
  (\id_A\oti\tilde d_A)\cir(\Delta\oti\id_\Av)\cir b_A= (\dimR(A)/\gam)\,\eta\,,
  \\{}\\[-.5em]
  d_A\cir(\id_\Av\oti m)\cir(\tilde b_A\oti\id_A) = (\dimL(A)/\gam)\,\eps \,,
  \\{}\\[-.5em]
  \tilde d_A\cir(m\oti\id_\Av)\cir(\id_A\oti b_A) = (\dimR(A)/\gam)\,\eps \,.
  \labl88 \eear  \end{equation}

              \newpage\noindent
These identities imply, in turn, that
  \be  \trr(\Delta\cir m) = (\dimR(A))^2_{}/\gam  \end{equation}
etc. When \A\ is in addition \tcc, there are analogous identities such as
  \be \bearl
  (\tilde d_A\oti\id_A)\cir(\id_A\oti\tilde c)\cir(\Delta\oti\id_\Av)
  \cir b_A = (\dimR(A)/\gam)\,\eta\,.  \labl89 \eear  \end{equation}

\vskip.85em

We have been careful to distinguish between $X$ and $X\Vee$ as objects of 
\calc\ even when $X$ is self-dual. It may seem desirable to identify any 
self-dual $X$ with $X\Vee$. But as the following considerations show, such
an identification is in general not possible. We will need the following

\medskip

\definition
An {\sl \absi\/} object $X$ of an abelian category is an object for which
  \be  \Hom(X,X) = k\,\id_X  \end{equation}
as a $k$-bimodule. 

\medskip

Further, for any self-dual object $X$ of a sovereign tensor category and any
isomorphism $\phi\iN\Hom(X,X\Vee)$ we define the endomorphism 
  \be  \Nu_X := (d_X\oti\id_{X\Vee}) \circ               
  (\id_{X\Vee}\oti\phi^{-1}\oti\phi) \circ (\id_{X\Vee}\oti\tilde b_X) 
  \;\in \Hom(X\Vee\!{,}\,X\Vee) \,. \labl nu  \end{equation}
This is in fact an automorphism, with inverse
  \be  \Nu_X^{-1} = (\tilde d_X\oti\id_{X\Vee}) \circ
  (\phi^{-1}\oti\phi\oti\id_{X\Vee}) \circ (\id_{X\Vee}\oti b_X) \,.
  \end{equation}
The morphisms $\Nu_X\,{\equiv}\,\Nu_X^{(\phi)}$ and $\Nu_X^{-1}$ may depend on
the choice of $\phi$. They can be regarded as functions on the $k$-projective 
space $P\Hom(X,X\Vee)$, where $k$ acts from the left. As a consequence, for 
\absi\ $X$ the morphism $\Nu_X$ does {\sl not\/} depend on 
the choice of $\phi$. This allows the

\definition
(i)\ The (right) {\sl Frobenius\hy Schur indicator\/} -- or {\sl 
FS indicator\/}, for short -- $\nu_X$ of a self-dual \absi\ 
object $X$ of a sovereign \tc\ is the invertible scalar appearing in
  \be  \Nu_X =: \nu_X\, \id_{X\Vee} \,,  \labl fs  \end{equation}
where $\Nu_X$ is the map \erf nu.
\\[.2em]
(ii)\ A self-dual \absi\ object $X$ of a sovereign \tc\ 
is called {\sl real\/} iff $\nu_X\eq1$.

\remark
1.\ Morally, the scalar $\nu_X$ measures the difference between
the duality morphisms $b_X$ and $(\phi^{-1}\Oti\,\phi)\cir\tilde b_X$. That
it can be different from 1 shows that even up to all possible isomorphisms 
left and right duality can be different for self-dual objects.
\\[.19em]
2.\ In the case of tensor categories constructed from $q$-deformed
universal enveloping algebras of \findim\ simple \lie s or from untwisted
affine \lie s, the FS indicator distinguishes orthogonal ($\nu\eq1$) from
symplectic ($\nu\eq{-}1$) simple objects.
\\[.19em]
3.\ In addition to the right FS indicator there is also a {\sl left
FS indicator\/} $\tilde\nu_X$ of \absi\ objects, defined by using the
isomorphism
  \be  \tilde\Nu_X^{(\phi)} := (\id_{X\Vee}\oti\tilde d_X) \circ
  (\phi\oti\phi^{-1}\oti\id_{X\Vee}) \circ (b_X\oti\id_{X\Vee})
  \labl nt  \end{equation}
of $X\Vee$ instead of $\Nu_X^{(\phi)}$.
\\[.19em]
4.\ The present definition of FS indicator generalizes the one of \cite{fffs3}. 
This follows from the fact that in a modular tensor category the left duality 
can be expressed 
as $\tilde b_X\eq(\theta_{X\Vee}\Oti\id_X)\cir c_{X,X\Vee}\cir b_X$
in terms of the right duality, the braiding $c$ and the twist $\theta$.
Compare also section 4 of \cite{fgsv} for the case of *-categories.

\smallskip

\proposition
The right Frobenius\hy Schur indicator satisfies
  \be  \dim(X)\, \nu_X^2 = \dim(X) \,. \labl2n  \end{equation}
In particular, $\nu_X^2\eq1$ if $\dim(X)$ is invertible.
\proof\dd
Multiplying $\Nu_X$ by $\nu_X$ and 
composing with $\tilde d_X$ from the left and with $b_X$ from the
right, one obtains directly $\nu_X^2\dimr(X)$. On the other hand, 
after replacing $\nu_X$ times the morphism $\id_{X\Vee}$ in the 
middle part of the definition \erf nu by $\Nu_X$, one pair of $\phi$ 
and $\phi^{-1}$ are concatenated and hence cancel; then after successive use 
of the duality axioms \erf3d and \erf2d one can use $\phi\cir\phi^{-1}
\eq\id_X$ once again so as to arrive at $\diml(X)$. More explicitly,
  \be \bearl  \nu_X^2\, \dimr(X) 
  \\{}\\[-.9em]\hsp{1.3}
  = \tilde d_X \circ (\id_X\oti d_X\oti\id_{X\Vee}) \circ
     (\id_X\oti\id_{X\Vee}\oti\phi^{-1}\oti\phi)
  \\{}\\[-.9em]\hsp{1.3}\quad
     \circ (\id_X\oti\id_{X\Vee}\oti d_X\oti\id_{X\Vee}\oti\id_X)
     \circ (\id_X\oti\id_{X\Vee}\oti\id_{X\Vee}\oti\phi^{-1}
     \oti\phi\oti\id_X)
  \\{}\\[-.9em]\hsp{1.3}\quad
     \circ (\id_X\oti\id_{X\Vee}\oti
     \id_{X\Vee}\oti\tilde b_X\oti\id_X) \circ (b_X\oti\tilde b_X)
  \\{}\\[-.9em]\hsp{1.3}
  = \tilde d_X \circ (\id_X\oti d_X\oti\phi) \circ \Llb
     \id_X\oti\id_{X\Vee}\oti (\phi^{-1}\circ\phi) \oti\id_X \Lrb
  \\{}\\[-.9em]\hsp{1.3}\quad
     \circ (\id_X\oti\id_{X\Vee}\oti d_X\oti\id_X\oti\id_X) \circ
     (\id_X\oti
     \id_{X\Vee}\oti\id_{X\Vee}\oti\phi^{-1}\oti\id_X\oti\id_X)
  \\{}\\[-.9em]\hsp{1.3}\quad
     \circ (\id_X\oti\id_{X\Vee}\oti\id_{X\Vee}\oti\tilde b_X\oti
     \id_X) \circ (b_X\oti\tilde b_X)
  \\{}\\[-.9em]\hsp{1.3}
  = \tilde d_X \circ \Llb \llb (\id_X\oti d_X)\circ(b_X\oti\id_X) \lrb
     \oti \phi \Lrb \circ (d_X\oti\id_X\oti\id_X)
  \\{}\\[-.9em]\hsp{1.3}\quad
     \circ (\id_{X\Vee}\oti\phi^{-1}\oti\id_X\oti\id_X)
     \circ (\id_{X\Vee}\oti\tilde b_X\oti\id_X) \circ \tilde b_X
  \\{}\\[-.9em]\hsp{1.3}
  = \tilde d_X \circ ((d_X\oti\id_X\oti\id_{X\Vee}) \circ
     (\id_{X\Vee}\oti\phi^{-1}\oti\id_X\oti\phi) \circ
     (\id_{X\Vee}\oti\tilde b_X\oti\id_X) \circ \tilde b_X
  \\{}\\[-.9em]\hsp{1.3}
  = d_X \circ (\id_{X\Vee}\oti\phi^{-1}) \circ \Llb \id_{X\Vee}
     \oti \llb (\id_{X\Vee}\oti\tilde d_X)\circ(\tilde b_X\oti
     \id_{X\Vee}) \lrb \Lrb \circ (\id_{X\Vee}\oti\phi) \circ \tilde b_X
  \\{}\\[-.9em]\hsp{1.3}
  = d_X \circ (\id_{X\Vee}\oti\phi^{-1}) \circ (\id_{X\Vee}\oti\phi)
     \circ \tilde b_X
  = d_X \circ \tilde b_X
  = \diml(X) \,.
  \\{}\\[-1.9em]
  \eear \end{equation}

\vskip.3em\noindent
The statement thus follows by equality of the left and right
dimensions (see formula \erf dd) of $X$.
\hfill\checkmark

\remark \labl35  
1.\ Analogously one obtains
  \be  \dim(X)\, \tilde\nu_X^2 = \dim(X)  \labl2m  \end{equation}
for the left FS indicator.
\\
2.\ In many cases of interest, in particular for every object of a 
modular category over $\complex$\,, the dimension is invertible.
Then we simply have $\nu_X^2\eq1\eq\tilde\nu_X^2$.
\\
3.\ A generalization of \erf2n to reducible objects $X$ is given by
  \be  \trl(\Nu_X^{(\phi)}) = \tilde d_X \circ (\phi^{-1}\oti\phi)
  \circ \tilde b_X = \trr((\Nu_X^{(\phi)})^{-1}_{}) \,.
  \labl3n  \end{equation}
Indeed, this relation reduces to \erf2n when $\Nu_X^{(\phi)}$ is a multiple
$\nu_X$ of $\id_{X\Vee}$.
\\
4.\ In particular, when \A\ is a Frobenius \alg\ in a sovereign category,
we can evaluate formula \erf3n 
using for $\phi$ the morphism $\tilde\Phi$ defined in \erf ft. We have
  \be  \Nu_{\!A}^{(\tilde\Phi)} = \tilde\Phi\cir\Phi^{-1}_{}
  \labl tp  \end{equation}
   and
  \be  \tilde\Nu_{\!A}^{(\Phi)} = \Phi\cir\tilde\Phi^{-1}_{}
  = (\Nu_{\!A}^{(\tilde\Phi)})^{-1}_{} \,.  \labl pt  \end{equation} 
We thus obtain
$\trl(\Nu_{\!A}^{(\tilde\Phi)})
                                \eq\trr(\tilde\Nu_{\!A}^{(\Phi)})
\eq\eps\cir m\cir\Delta\cir\eta$, or what is the same,
  \be  \trl(\tilde\Phi\cir\Phi^{-1}_{}) = \trr(\Phi\cir \tilde
  \Phi^{-1}_{}) = \eps\circ m\circ\Delta\circ\eta \,.  \end{equation}
When \A\ is special Frobenius, this reduces to
  \be  \trl(\tilde\Phi\cir\Phi^{-1}_{}) = \trr(\Phi\cir \tilde
  \Phi^{-1}_{}) = \bet\gam \,.  \labl bg  \end{equation}
We also note the identity
  \be  d_A \circ (\tilde\Nu_{\!A}^{(\Phi)} \oti (\eta\cir\eps))
  \circ \tilde b_A = \eps \circ m \circ (\id_A \oti (\eta\cir\eps))
  \circ \Delta \circ \eta = \eps \circ \eta = \gam \,.  \labl8x
  \end{equation}

             \vfill\noindent
5.\ Whenever one has $\Nu_{\!A}\eq\id_{A\Vee}$, 
the identity \erf bg implies
  \be  \bet\gam = \dimL(A) \,.  \labl a1  \end{equation}
Also note that the notion of reality should be extendible to apply to
non-\absi\ self-dual objects. In general one should {\em not\/} expect that
the existence of an isomorphism $\phi$ with $\Nu_{\!X}^{(\phi)}
\eq\id_{X\Vee}$ is necessary for $X$ to be real.
(But when \A\ is in addition \hapl\ in the sense of definition
\ref{11} below, then more can be said; see the formulae \erf52 and \erf a2.)

             \vfill\noindent
6.\ When applied to the tensor unit $\one$ -- which is a special Frobenius 
\alg, see remark \ref{st}.3 -- these results imply that
the dimension as well as the FS indicator of $\one$ are equal to one,
  \be  \dim(\one) = 1 \qquad{\rm and}\qquad \nu_\smallone = 1 \,.
  \labl=1  \end{equation}
Indeed, the traces \erf bg are then equal to 1. But the morphisms 
$\Nu_{\!\smallone}^{(\tilde\Phi)}$ and $\tilde\Nu_{\!\smallone}^{(\Phi)}$
are multiples of $\id_{\smallone\Vee}$. It follows that
$\nu_\smallone\eq(\dim(\one))^{-1}$, which when combined with the result 
\erf2n implies $(\dim(\one))^2\eq1$. On the other hand, applying
the general identity $\diml(X\Oti Y)\eq\diml(X)\,\diml(Y)$ (which follows
by taking the trace of $\id_{X\otimes Y}\eq\id_X\Oti\id_Y$) to the case
$X\eq\one$ and $Y\eq\one\Vee$ one gets $(\dim(\one))^2\eq\dim(\one)$. 
By combining the two results, we end up with \erf=1.

                    \vfill\vfill
\section{Modules}

\definition
(i)\ For \A\ an \alg\ in \calc, a {\sl left $A$-module\/} in a tensor 
category $\calc$ is a pair
$(N,\R)$ with $N\iN\obj(\calc)$ and $\R\iN\Hom(A\Oti N,N)$ satisfying
  \be  \R\circ(m\oti\id_N) = \R \circ (\id_A\oti\R)
  \qquad\mbox{and}\qquad \R\circ(\eta\oti\id_N) = \id_N \,.  \labl1m
  \end{equation}
(ii)\ For every $X\iN\obj(\calc)$, the {\sl induced\/} (left) {\sl module\/}
$\inda(X)$ is the pair 
  \be  (A\Oti X,m\Oti\id_X) \,.  \end{equation}
(iii)\ For \A\ an \alg\ in \calc, the {\sl category \calca\ of $A$-modules\/}
in \calc\ is the category whose objects are the $A$-modules in \calc\
and whose morphisms are
  \be  \Hom_A(N,M) := \{ f\iN \Hom(N,M) \,|\,
  f \cir\R_N\eq\R_M\cir(\id_A\Oti f) \} \,.  \labl ha  \end{equation}
(iv) An \A-module $(N,\R_N)$ is called {\sl \absi\/} iff 
$\Hom_A(N,N)\eq k\,\id_N$.  

                   \newpage 

\remark \labl42
1.\ Right modules and bimodules, as well as left and right comodules, are 
defined analogously. Below the term module will always refer to left modules.
In the second of the relations \erf1m it is used that $\one\oti N$
is the same object as $N$. 
\\
2.\ That $\inda(X)$ is a (left) \A-module is a direct consequence of the
associativity \erf as and of the unit property \erf un.
Analogously one can define an induced right module as
$\indar(X)\,{:=}\,(X\Oti A,\id_X\Oti m)$, which is a right \A-module.
\\
\A\ is both a left and a right module over itself; indeed we have the 
equalities
  \be  \inda(\one) = A \qquad{\rm and}\qquad \indar(\one) = A
  \end{equation}
of left and right \A-modules, \resp.
The module properties are nothing but the defining axioms of the \alg\ 
structure. Moreover, the left and right actions of \A\ on itself 
commute. Indeed, this property is equivalent to the associativity of $m$ 
(it is worth noting that commutativity of $m$ is {\sl not\/} required).
\\
3.\ Not every object $X$ of $\calc$ needs to underly some $A$-module,
and an object of $\calc$ can be an $A$-module in several inequivalent ways.
\\
4.\ We will often suppress the morphism $\R$ and just write $N$ for a
(left or right) $A$-module $(N,\R)$. When necessary, we distinguish the 
module $N$ from the underlying object of \calc\ by denoting the latter as
$\Dot N$.
\\
5.\ When $N$ is a submodule of $M$ (i.e.\ a subobject in \calca),
then $\Dot N$ is a subobject of $\Dot M$ in \calc.
\\
6.\ The product of \A\ is a morphism in \calca,
  \be  m \in \Hom_A(A\Oti A,A)\,.  \labl ma  \end{equation}
The coproduct, on the other hand, is a morphism in \calca, i.e.
  \be  \Delta \in \Hom_A(A,A\Oti A)\,,  \labl da \end{equation}
only if \A\ is Frobenius. 
Indeed, \erf ma is equivalent to the associativity \erf as, while \erf da 
is equivalent to the second half of the Fro\-be\-ni\-us axiom \erf1f
(and the full Fro\-be\-ni\-us axiom to the assertion that $\Delta$ is a
morphism of $A$ and $A\Oti A$ as \A-bimodules).
\\
7. \calca\ is indeed a category: Associativity of the concatenation 
of morphisms in \calca\ follows directly from associativity in \calc,
and for all modules $(N,\R)$ the morphism $\id_N$
is in the subspace $\Hom_A(N,N)$ of $\Hom(N,N)$. 
\\
8.\ The subsets \erf ha of morphisms of \calc\ are certainly the ones
of most direct interest in the context of our construction. However,
from the treatment of CFT \bc s in \cite{fuSc112} it is apparent that 
for various purposes one will also be interested in `twisted versions'
of the concept, i.e.\ with subsets
  \be  \Hom_A^{\{\varphi\}}(N,M) := \{ f\iN \Hom(N,M) \,|\,
  f \cir\R_N\eq\R_M\cir(\varphi\Oti f) \}  \end{equation} 
for suitable endomorphisms $\varphi\iN\Hom(A,A)$. Particular examples of 
such morphisms have been studied in \cite{kirI14}. At the representation 
theoretic level, such maps correspond to twisted intertwiners. They play 
an important role in the treatment of orbifold theories (see e.g.\ 
\cite{bifs}); their traces, so-called twining characters \cite{fusS3,furs}
frequently show an interesting behavior under modular transformations.

\vskip.53em

To endow the collection of \A-modules in \calc\ with more structure, it can
be  useful to require \A\ to possess an additional property. This is 
formulated in

\definition \labl11
A {\sl \hapl\/} \alg\ \A\ in \calc\ is an \alg\ in \calc\ such that the set of
morphisms between $\one$ and \A\ is isomorphic to $k$ as a $k$-bimodule.

\vskip.2em

\noindent
Given this property, the morphisms between $\one$ and \A\ form in particular 
a free $k$-module of rank
  \be  \dim\, \Hom(\one,A) = 1 \,.  \labl1i  \end{equation}

\remark
1.\ When $X\iN\obj(\calc)$ is \absi, then the object $X\Oti X\Vee$, which 
according to lemma \ref{3r} is a special Frobenius algebra, is also \hapl.
\\
2.\ The special Frobenius \alg s\ ${\mathcal F}{\rm G}$ 
described in remark \ref{st}.4 are \hapl, too; this follows from the fact
that the trivial \rep\ of G is one-dimensional.

\vskip.3em

Without explicit mentioning, below we will often make the additional\,%
 \footnote{~In \cite{kios} \hapy\ of \A\ is included in the definition of
 the term `algebra'.}
assumption that the \alg\ \A\ of our interest is \hapl.
In the context of extensions of vertex algebras, this requirement is
entirely natural, since it encodes the uniqueness of the vacuum. (Concretely,
the vector space \aext\ possesses a distinguished element $v_\Omega$, the 
vacuum vector, and there is a unique irreducible \abar-module that contains 
this vector, namely \abar\ itself. The condition \erf1i says that \abar\ 
appears with multiplicity one in the decomposition of \aext\ as an 
\abar-module.) The following observation shows that
this assumption is also natural from a mathematical point of view.

\lemma
Every \hapl\ \calc-\alg\ \A\ is an \absi\ left- and 
right module over itself. In particular,
  \be  \Hom_A(A,A) = k\, \id_A \,.  \labl1j  \end{equation}

\proof\dd
Because of \erf1i every $\varphi\iN\Hom_A(A,A)$ satisfies
$\varphi\cir\eta\eq\xi\,\eta$ for some $\xi\iN k$. Thus we can
use the unit property of $\eta$ and the module property to see that 
  \be  \bearll  \varphi \!\!\!
  &= \varphi \circ m \circ (\id_A\oti\eta)
   = m \circ (\id_A\oti\varphi) \circ (\id_A\oti\eta)
  \\{}\\[-.7em]
  &= m \circ (\id_A\oti (\varphi\circ\eta))
   = m \circ (\id_A\oti \xi\,\eta) = \xi\,\id_A \,.
  \eear  \end{equation}
{}\mbox{$\ $}\\[-1.5em]\phantom.
\hfill\checkmark

\vskip.3em
\remark
1.\ The first assertion is also a direct consequence of the reciprocity
relation \erf4g that will be discussed below.
\\
2.\
It follows in particular that $m\cir\Delta$ is a scalar multiple of $\id_A$.
Thus as soon as the restriction \erf1i is made, the non-trivial statement 
in the special property \erf3f is that the scalar $\bet$ is invertible.  

\vskip.53em
To proceed, let us introduce an operation $E$ on morphisms by
  \be  E(f) \equiv E_{X,N}(f)
  := \R_N \circ (\id_A\oti f) \;\in \Hom(A\Oti X,\Dot N)  \labl ee
  \end{equation}
for all $X\iN\obj(\calc)$, $N\iN\obj(\calca)$ and $f\iN\Hom(X,\Dot N)$.
$E$ is used in establishing the

\proposition
For \A\ an \alg\ in \calc, $X$ an object of\, \calc\ and $N$ an \A-module 
in \calc, there is a canonical isomorphism
  \be  E:\quad  \Hom(X,\Dot N) \,\stackrel\cong\to\, \Hom_A (A\Oti X,N) \,.
  \labl1k  \end{equation}

\proof\dd
We give $A\Oti X$ the structure of an induced \A-module. Then we
relate elements of $\Hom(X,N)$ and $\Hom_A(A\Oti X,N)$ by the
operation $E$ defined in formula \erf ee and by setting
  \be  R(g) := g \circ (\eta\oti\id_X)  \;\in \Hom(X,N)  \labl46
  \end{equation}
for $g\iN\Hom_A(A\Oti X,N)$. Next we note that
the morphism $E(f)$ is in fact in the set $\Hom_A(A\Oti X,N)$, because
owing to the representation property of $\R_N$ we have
  \be \bearll
  E(f) \circ \R_{A\oti X} \!\!\!
  &= \R_N \circ(\id_A\oti f) \circ (m\oti \id_X) = \R_N \circ (m\oti f)
  \\{}\\[-.9em]
  &= \R_N \circ (m\oti \id_N) \circ (\id_A\oti\id_A\oti f)
  \\{}\\[-.9em]
  &= \R_N \circ (\id_A\oti\R_N) \circ (\id_A\oti\id_A\oti f) 
   = \R_N \circ (\id_A\oti E(f)) \,.
  \eear  \end{equation}
Furthermore, the two operations $E$ and $R$ are each other's inverse. 
Namely, using the functoriality of the tensor product
and the representation property of $\R_N$ we find
  \be \bearll  R(E(f))  \!\!\!
  & 
   = \R_N \circ (\id_A\oti f) \circ (\eta\oti\id_X)
   = \R_N \circ (\eta\oti f) 
  \\{}\\[-.9em]
  &= \R_N \circ(\eta\oti\id_N) \circ (\id_{\smallone}\oti f)
   = \id_N \circ (\id_{\smallone}\oti f) = f
  \,, \eear  \end{equation}
while from the fact that $g$ intertwines the induced \rep\ with $N$,
the explicit form of the induced action and the unit property of $\eta$
it follows that
  \be \bearll  E(R(g)) \!\!\!
  & 
   = \R_N \circ (\id_A \oti g) \circ (\id_A\oti\eta\oti\id_X) 
  \\{}\\[-.9em]
  &= g \circ (m\oti\id_X) \circ (\id_A\oti\eta\oti\id_X) 
   = g \circ (\id_A\oti \id_X) = g
  \,. \eear  \end{equation}
{}\mbox{$\ $}\\[-1.4em]\phantom.\hfill\checkmark
\vskip.3em

Note that as a special case of the result \erf1k, it follows that
when we do not make the assumption that the \alg\ \A\ is \hapl,
then the relation \erf1j gets relaxed to
  \be  \Hom_A(A,A) \equiv \Hom_A(\inda(\one),A) \cong \Hom(\one,A)  \,.
  \end{equation}
Also, the result motivates the following

\definition
(i)\ The {\sl induction functor\/} $\F{:}\ \calc\,{\to}\,\calca$ is
defined by
  \be  \F(X) := \inda(X) \qquad{\rm for}\quad X\iN\obj(\calc)
  \end{equation}
and
  \be  \F(f) := E_{X,A\otimes Y}(\eta\oti f) = \id_A \oti f
  \qquad{\rm for}\quad f\iN\Hom(X,Y) \,.  \end{equation}
(ii)\ The {\sl restriction functor\/} $\G{:}\ \calca\,{\to}\,\calc$ is
the forgetful functor, i.e.
  \be  \G(N) := \Dot N \qquad{\rm for}\quad N\equiv(\Dot N,\R_N)
  \iN\obj(\calca)  \end{equation}
and for $g\iN\Hom_A(M,N)$, $\G(g)$ is the same morphism considered as
an element of $\Hom(\Dot M,\Dot N)\,{\supseteq}\,\Hom_A(M,N)$.

\lemma
$\F$ and $\G$ are functors. When \calc\ is autonomous
and \A\ is Frobenius, then these functors are exact.

\proof\dd
For every $X\iN\obj(\calc)$ we have $\F(\id_X)\eq\id_A\oti\id_X\eq
\id_{A\otimes X}\eq\id_{\F(X)}$ and, for $f\iN\Hom(X,Z)$ and $g\iN\Hom(Z,Y)$,
  \be  \F(g\cir f) = \id_A\oti(g\cir f) = (\id_A\oti g) \circ
  (\id_A\oti f) = \F(g) \circ \F(f) \,.  \end{equation}
Thus $\F$ is a functor.  When \A\ is Frobenius,
then it makes sense to ask about ex\-act\-ness because, as will be shown in
Proposition \ref{ba} below, in that case the module category \calca\ is 
abelian. Exactness of $\F$ then follows immediately from bi-exactness
of the tensor product functor $\otimes$ in an autonomous \tc.
\\  
For $\G$ the statements are trivial. 
\hfill\checkmark

\vskip.53em

\proposition
$\F$ is a right adjoint functor of \,$\G$.

\proof\dd
According to proposition XI.1.8.\ of \cite{KAss} we must show that the 
bijections
  \be  R: \quad \Hom_A(\F(X),M) \to \Hom(X,\G(M)) \labl4g  \end{equation}
that are provided by formula \erf46 make the diagram
  \be \begin{array}{ccc}
  \Hom_A(\F(Y), M) &\ \ \stackrel{\scriptstyle R}{\longrightarrow}
  \ \ & \Hom(Y,\G(M)) \\{}\\[-.5em]
  \quad \downarrow {\scriptstyle \circ\,\F(g)} & &
  \downarrow {\scriptstyle \circ\, g} \\{}\\[-.5em]
  \Hom_A(\F(X), M) &\ \ \stackrel{\scriptstyle R}{\longrightarrow}
  \ \ & \Hom(X,\G(M)) \\{}\\[-.5em]
  \downarrow {\scriptstyle f\circ} & &
  \quad \downarrow {\scriptstyle \G(f)\,\circ} \\{}\\[-.5em]
  \Hom_A(\F(X),N) &\ \ \stackrel{\scriptstyle R}{\longrightarrow}
  \ \ & \Hom(X,\G(N)) \\{}\\[-.5em]
  \eear  \end{equation}
commutative for every $X,Y\iN\obj(\calc)$, $M,N\iN\obj(\calca)$, 
$f\iN\Hom_A(M,N)$ and $g\iN$\linebreak[0]$\Hom(X,Y)$. By direct calculation 
one checks that this is indeed the case: The upper square commutes because 
for each $h\iN\Hom_A(A\oti Z,M)$ we have the equality
  \be  R(h)\circ g
  = h \circ (\eta \oti \id_Z) \circ g 
  = h\circ (\id_A\oti g) \circ (\eta\oti\id_X) 
  = R(h\cir\F(g))  \end{equation}
as morphisms in $\Hom(X,M)$. Similarly,
  \be  \G(f) \circ R(h') = f \circ h' \circ (\eta\oti\id_X)
  = R(f\cir h')  \end{equation}
for every $h'\iN\Hom(A\oti X,M)$, so that the lower square commutes as well.
\hfill\checkmark

\vskip.63em

When the \alg\ is Frobenius, then in addition to the reciprocity relation
\erf4g there is also the following mirror version.

\proposition \labl*0
When \A\ is Frobenius, then for every $M\iN\obj(\calca)$ and every
$X\iN\obj(\calc)$ there is a natural bijection
  \be  \tilde E:\quad \Hom_A(M,\F(X)) \cong \Hom(\Dot M,X)
  \labl55  \end{equation}
of $k$-bimodules.

\proof\dd
We set
  \be  \tg f := (\eps\oti\id_X)\circ f
  \;\in\Hom(\Dot M,X)  \end{equation}
for any $f\iN\Hom_A(M,\F(X))\eq\Hom_A(M,\inda(X))$ and 
  \be  \tf g := (\id_A\oti g) \circ (\id_A\oti\R_M) \circ ((\Delta\cir\eta)
  \oti\id_{\Dot M})  \;\in\Hom(\Dot M,A\Oti X)  \end{equation}
for any $g\iN\Hom(\G(M),X)\eq\Hom(\Dot M,X)$. One sees that
the morphisms $\tf g$ are in fact in $\Hom_A(M,\inda(X))$ by using 
successively the \rep\ property of $\R_M$, the Frobenius property, the unit
property of $\eta$, again the Frobenius property, and the definition of the 
induced module $\inda(X)$. Further, by using the counit property and then the 
\rep\ property of $R_M$, it follows that $\tg{\tf g}\eq g$ for all 
$g\iN\Hom(\Dot M,X)$. And that $\tf{\tg f}\eq f$ for all $f\iN$\linebreak[0]
$\Hom_A(M,\inda(X))$ holds as well follows because, by the definition of the
$\Hom_A$ sets and of $\inda(X)$, one has
  \be  \tf{\tg f} = (\id_A\oti(\eps\cir m)\oti\id_{\Dot M})
  \circ ((\Delta\cir\eta)\oti f) \,;  \end{equation}
this reduces to $f$ by the Frobenius axiom and the unit and counit
properties.
\hfill\checkmark

\vskip.53em

\lemma \labl*1
When \A\ is a Frobenius \alg, then for any $M,N\iN\obj(\calca)$, $X,Y\iN\obj
(\calc)$ and $f\iN\Hom(\Dot M,X)$ the following statements hold.\\
First, for all $h\iN\Hom(Y,\Dot M)$, the morphisms $\tf f\iN\Hom_A(M,\inda(X))$
and $E(h)$\linebreak[0]$\iN\Hom_A(\inda(Y),N)$ obey
  \be  \tf f \circ E(h) = (\id_A \oti (f\cir\R_M)) \circ (\Delta\cir h) \,.
  \labl fh  \end{equation}
Second, when \A\ is also special, then for all morphisms $g\iN\Hom(X,
\Dot N)$ that satisfy $g\cir f\iN\Hom_A(M,N)$ one has
  \be  E(g) \circ \tf f = \bet\, g \circ f \,.  \labl gf  \end{equation}

\proof\dd 
The first relation follows by using the \rep\ property of $\R_M$ and the
Frobenius axiom. As to the second, we have by definition
  \be  E(g) \circ \tf f = \R_N \circ (\id_A \oti (g\cir f\cir \R_M))
  \circ ((\Delta\cir\eta) \oti \id_{\Dot M}) \,.   \end{equation}
Using the fact that $g\cir f$ is a morphism in \calca, this can be
rewritten as
  \be  E(g) \circ \tf f = g \circ f \circ \R_M \circ (\id_A \oti \R_M)
  \circ ((\Delta\cir\eta) \oti \id_{\Dot M}) \,.   \end{equation}
The result then follows by using the \rep\ property of $\R_M$, the
special property of \A\ and again the \rep\ property of $\R_M$.
\hfill\checkmark

\vskip.3em

\lemma \labl*2
For every $M,N\iN\obj(\calca)$, $X\iN\obj(\calc)$ and for all
morphisms $f\iN\Hom_A(M,N)$ the following holds. First, for all
$h\iN\Hom(X,\Dot M)$ one has
  \be  E(f\cir h) = f \circ E(h) \,. \labl22  \end{equation}
Second, when \A\ is a Frobenius \alg, then
  \be  \tf{g\cir f} = \tf g \circ f  \labl21  \end{equation}
for all $g\iN\Hom(\Dot N,X)$.

\proof\dd
Using the fact that $f$ is a morphism in \calca, this follows immediately
from the definitions:
  \be  \bearll 
  E(f\cir h) \!\!\!\!
  &= \R_N \circ (\id_A\oti(f\cir h))
   = f \circ \R_M \circ (\id_A\oti h) = f \circ E(h)
  \,,\\{}\\[-.6em]
  \tf{g\cir f} \!\!\!\!
  &= (\id_A\oti(g\cir f\cir \R_M)) \circ ((\Delta\cir\eta)\oti\id_{\Dot M}) 
  \\{}\\[-.8em]
  &= (\id_A\oti(g\cir\R_N)) \circ ((\Delta\cir\eta)\oti f)
  \\{}\\[-.8em]
  &= (\id_A\oti(g\cir\R_N)) \circ ((\Delta\cir\eta)\oti\id_{\Dot N}) \circ f
   = \tf g \circ f 
  \,. \eear  \end{equation}
(Recall that the Frobenius property is used in the definition of $\tilde E$.)
\hfill\checkmark

\vskip.3em

\lemma \labl*3
{\rm(i)}~Every monomorphism of \calca\ is also a monomorphism of \calc.
\\
{\rm(ii)}~When \A\ is a Frobenius \alg,
then every epimorphism of \calca\ is also an epimorphism of \calc.

\proof\dd
(i) Let $f\iN\Hom_A(M,N)$ be a monomorphism of \calca, and assume that
$g,h\iN\Hom(X,\Dot M)$ for some $X\iN\obj(\calc)$ obey $f\cir g\eq f\cir h$
(where $f$ is regarded as a morphism of \calc). Then also
the corresponding induced morphisms are equal, $E(f\cir g)\eq E(f\cir h)$,
and hence $f\cir E(g)\eq f\cir E(h)$ by formula \erf22. Since $f$ is a
monomorphism of \calca, it follows that $E(g)\eq E(h)$, which in turn implies 
that $g\eq h$ because induction is bijective on morphisms (see \erf1k). 
Since this result holds for all $g$ and $h$, $f$ is a monomorphism of \calc.
\\
(ii) Similarly, let $f\iN\Hom_A(M,N)$ be an epimorphism of \calca, and 
assume that $g,h\iN\Hom(\Dot N,Y)$ for some $Y\iN\obj(\calc)$ obey 
$g\cir f\eq h\cir f$. Then with the help of formula \erf21 it follows that
$\tf g\cir f\eq\tf{g\cir f}\eq\tf{h\cir f}\eq\tf h\cir f$. Then $\tf g\eq\tf h$
by the epimorphism property of $f$, and hence $g\eq h$ by \erf55.
\hfill\checkmark

\lemma  \labl99
When \A\ is special Frobenius, then every \A-module is a submodule of an
induced module.

\proof\dd
Using the structure of $A\oti A$ as a right $A$-module and the definition
of $\Ota$ as an image, it is easy to see that for every $M\iN\obj(\calca)$ 
one has an isomorphism
between $(A\Oti A)\,\Ota M$ and $A\oti(A\Ota M)$, and thereby also between
$(A\Oti A)\,\Ota M$ and $A\oti M\eq \F(\Dot M)$, as left $A$-modules. Moreover, 
since $A$ is special Frobenius, \A\ is a sub-bimodule of the \A-bimodule 
$A\Oti A$. As a consequence $M\,{\cong}\,A\,\Ota M$ is a submodule of 
$(A\Oti A)\,\Ota M$ and hence, by the isomorphism just mentioned, a submodule 
of the induced module $\F(\Dot M)$.
\hfill\checkmark

\section{Some properties of the module category \calca}

The category \calca\ would not suit our purposes if it had no additional 
structure. In the sequel we demonstrate that when the \calc-\alg\ \A\ is 
sufficiently well behaved, then the module category \calca\ inherits 
various structural properties from \calc.

\subsection{Abelianness}
{}\mbox{$\ $}\\[-1.0em]
\proposition \labl ba
When \A\ is a special Frobenius \alg\ in \calc,
then the module category \calca\ is abelian.

\proof\dd 
(i) By construction \calca\ is $k$-additive. Using that the morphisms
$\Hom(A$\linebreak[0]$\Oti X_0,X_0)$ involving the zero object $X_0$ 
of \calc\ consist of a single element, one sees that $X_0$ is an \A-module, 
which provides a zero object of \calca.
\\
(ii) Next, let $f\iN\Hom_A(M,N)$. Then $f$ is in
$\Hom(\Dot M,\Dot N)$ and hence possesses a kernel $(\Dot K,\jmath)$ as a
morphism of \calc, i.e.\ there exist $\Dot K\iN\obj(\calc)$ and
$\jmath\iN\Hom(\Dot K,\Dot M)$ which obey $f\cir\jmath\eq0$ as well as the
universal property that for every $g\iN\Hom(\Dot M',\Dot M)$ with $f\cir g\eq0$
there is a unique $\Dot h\iN\Hom(\Dot M',\Dot K)$ with $\jmath\cir\Dot h\eq g$.
Since $f$ is a module morphism, the morphism
  \be
  \tilde\rho_K := \rho_M \cir (\id_A \oti \jmath) ~\in\Hom(A\oti\Dot K,\Dot M)
  \eee
satisfies
  \be
  f \cir \tilde\rho_K = \rho_N \cir (\id_A \oti (f\cir\jmath)) = 0 \,.
  \eee
By the universal property of $(\Dot K,\jmath)$ there is thus a unique
$\rho_K \iN \Hom(A\oti\Dot K,\Dot K)$ such that
  \be
  \jmath \circ \rho_K = \tilde\rho_K \,.
  \label{roK}\eee
Using the module properties of $M$ and applying the universal property
of $(\Dot K,\jmath)$ to the morphisms $\rho_K \cir (\id_A\oti\rho_K)$
and $\rho_K \cir (m\oti \id_A)$ as well as to $\rho_K \cir
(\eta \oti \id_{\dot K})$ and $\rho_K \cir (\id_{\dot K} \oti \eta)$,
one verifies that
  \be
  K := (\dot K,\rho_K)
  \eee
defines a left $A$-module structure on the object $\Dot K$. Further, with
the help of (\ref{roK}) one checks that $\jmath$ is a module morphism,
$\jmath \iN \Hom_A(K,M)$.
\\
Now assume that $M'\iN \obj(\calca)$ and $\tilde g \iN \Hom_A(M',M)$ such
that $f \cir \tilde g \eq 0$. Since $\tilde g$ is an element of
$\Hom(\Dot M',\Dot M)$, there exists a unique
$g \iN \Hom(\Dot M',\Dot K)$ obeying $\jmath \cir g \eq \tilde g$.
Also, $f \cir \big( \rho_M \cir (\id_A \oti \tilde g) \big)
\eq \rho_N \cir \big( \id_A \oti (f \cir \tilde g) \big) \eq 0$, so
there is a unique $h \iN \Hom(A\oti\Dot M',\Dot K)$
such that $h \cir \jmath \eq \rho_M \cir (\id_A \oti \tilde g)$.
But this equality is satisfied
both by $h_1 \,{:=}\, g \cir \rho_{M'}$ and by $h_2 \,{:=}\, \rho_K \cir g$,
so $h_1 \eq h_2$, and hence $g \iN \Hom_A(M',K)$.
Thus $(K,\jmath)$ possesses the required universal property, and hence is a
kernel of $f$ in \calca.
\\
(iii) Analogously one shows that every morphism $f\iN\Hom_A(M,N)$ of 
\calca\ has a cokernel. It is provided by the module
$\inda(C)$ and the morphism $\tf p$, where $C\iN\obj(\calc)$ and 
$p\iN\Hom(\Dot N,C)$ are the object and projection that furnish the
cokernel of $f$ as a morphism of \calc. 
\\
(iv) Given a morphism $f\iN\Hom_A(M,N)$, 
we regard it as a morphism in $\Hom(\Dot M,\Dot N)$. Then by abelianness of
the category \calc\ there exists an $X\iN\obj(\calc)$ and a monomorphism
$i\iN\Hom(\Dot M,X)$ as well as an epimorphism $p\iN\Hom(X,\Dot N)$ such
that $f\eq i\cir p$. Now assume that there are $g,h\iN\Hom_A(M,\inda(X))$
such that $E(i)\cir g\eq E(i)\cir h$. Then because of proposition \ref{*0}
there also exist $\Dot g,\Dot h\iN\Hom(\Dot M,X)$ such that $g\eq\tf{\Dot g}$
and $h\eq\tf{\Dot h}$, and hence $E(i)\cir\tf{\Dot g}\eq E(i)\cir\tf{\Dot h}$.
By lemma \ref{*1}, this implies that $i\cir\Dot g\eq i\cir\Dot h$. Since
$i$ is a monomorphism in \calc, it follows that already $\Dot g\eq\Dot h$,
and hence also $g\eq h$. This shows that $E(i)$ is a monomorphism in \calca.
Analogously one can see that $\tf p$ is an epimorphism in \calca.
Finally, by lemma \ref{*1} the induced morphisms also
satisfy $E(i)\cir\tf p\eq\bet\,i\cir p$. Together these results tell us 
that in \calca\ every morphism $f$ can be written as the composition 
$f\eq j\cir q$ of a monomorphism $j$ and an epimorphism $q$.
\\
(v) Let $f\iN\Hom_A(M,N)$ be a monomorphism of \calca, $p\iN\Hom_A(N,C)$ 
the projection to the cokernel of $f$, and $\hilde g\iN\Hom_A(N',N)$ such
that $p\cir\hilde g\eq0$. {}From (iii) we know that there are 
$B\iN\obj(\calc)$ and $q\iN\Hom(\Dot N,B)$ such that $C\eq\inda(B)$ and
$p\eq\tf q$. Using \erf21 we thus have $\tf{q\cir\hilde g}\eq p\cir\hilde g
\eq0$, and hence $q\cir\hilde g\eq0$. Since $f$ is, by lemma \ref{*3}, 
a monomorphism even in \calc, abelianness of \calc\ then implies that
$f$ is the kernel of $p$ in \calc. As a consequence, there exists a unique
$g\iN\Hom(\Dot N',\Dot M)$ such that $f\cir g\eq\hilde g$. This morphism is
even a morphism of \calca. Indeed, using that $f$ and $\hilde g$ are
morphisms of \calca\ one finds
  \be  \bearll  f \circ \R_M \circ (\id_A \oti g) \!\!\!\!
  &= \R_N \circ (\id_A \oti (f\cir g))
  \\{}\\[-.8em]  
  &\equiv \R_N \circ (\id_A \oti \hilde g)
   = \hilde g \circ \R_{N'} \equiv f \circ g \circ \R_{N'} 
  \,, \eear  \end{equation}
which implies $\R_M\cir(\id_A\Oti g)\eq g\cir\R_{N'}$ because $f$ is 
a monomorphism. To summarize, for every $\hilde g\iN\Hom_A(N',N)$ with
$p\cir\hilde g\eq0$ there exists a unique $g\iN\Hom_A(N',M)$ obeying
$f\cir g\eq\hilde g$. Thus $f$ is the kernel of its cokernel.
\\
(vi) Analogously one shows that every epimorphism in \calca\ is the cokernel 
of its kernel.
\\
We have thus verified all properties that are needed for the module category
\calca\ to be abelian.
\hfill\checkmark

\vskip.3em
The following identity is trivial (it holds analogously in any abelian 
category), but we note it for later reference. Let $p\iN\Hom_A(N,C)$ be 
the cokernel projection for some morphism $f\iN\Hom_A(N',N)$, and let
$g^\flat\iN\Hom_A(M,N)$ and $\hilde h\iN\Hom_A(N,M')$ be arbitrary morphisms.
Introduce $g\,{:=}\,p\cir g^\flat$ and denote by $h\iN\Hom_A(C,M')$ the unique
morphism satisfying $h\cir p\eq\hilde h$. Then
  \be  h \circ g = h \circ p \circ g^\flat = \hilde h \circ g^\flat \,.
  \labl gh  \end{equation}

\vskip.4em

\subsection{Tensoriality}
{}\mbox{$\ $}\\[.4em] 
Our next aim is to obtain a module category that is again a \tc. 
We start by introducing the notion of a tensor product of modules. 
Just like in commutative \alg\ (or for arbitrary \alg s over $\complex$),
the idea is to equate the (right) action of \A\ on the first factor of a
tensor product of objects of \calc\ with the (left) action on the second
factor. Without making further assumptions on the modules we can only give 
the following

\definition
The {\sl tensor product\/} $M\ota N$ of a right \A-module $M$ and a 
left \A-module $N$ is the object
  \be  M\ota N := \coker\llb (\R_M \oti \id_N) - (\id_M \oti \R_N) \lrb
  \ \in\obj(\calc)  \labl0a  \,.  \end{equation} 
The corresponding projection is denoted by
  \be  p_{M,N} \in \Hom(M\Oti N,M\Ota N) \,.  \labl pc  \end{equation}

\vskip.53em

The \alg\ \A\ itself serves as a left and right unit for this tensor
product. This follows as a special case of the

\lemma \labl0l
For every left \A-module $M$ and every $X\iN\obj(\calc)$ one has an
isomorphism 
  \be  \indar(X)\ota M \,\cong\, X\oti\Dot M   \labl01  \end{equation}
as objects in \calc. Analogously, when $N$ is a right \A-module, then
  \be  N \ota \inda(X) \,\cong\, \Dot N\oti X  \labl02  \end{equation}
as objects in \calc.

\proof\dd
To establish the first isomorphism we must show that the object $X\Oti\Dot M$
is the cokernel of the morphism
  \be  f := \id_X\oti (m\oti\id_{\Dot M} - \id_A\oti\R_M) \
  \in\Hom(X\Oti A\Oti A\Oti\Dot M,X\Oti A\Oti\Dot M) \,,  \end{equation}
or more explicitly, that there is a morphism $p\iN\Hom(X\Oti A\Oti\Dot M,
X\Oti\Dot M)$ such that $p\cir f\eq0$ and such that for any 
$\hilde g\iN\Hom(X\Oti A\Oti\Dot M,Y)$ with $\hilde g\cir f\eq0$ there 
exists a unique $g\iN\Hom(X\Oti\Dot M,Y)$ such that $g\cir p\eq\hilde g$.
Clearly, the morphism $p\,{:=}\,\id_X\Oti\R_M$ satisfies the first 
requirement. Further, assume that for given $\hilde g$ there is a $g$
with the desired property. Then by the module property of $M$ we have
  \be  g = g \circ \llb \id_X\oti (\R_M\cir(\eta\Oti\id_{\Dot M})) \lrb
  = \hilde g \circ (\id_X\oti\eta\oti\id_{\Dot M}) \,.  \labl07
  \end{equation}
Thus in particular $g$ is unique. To show that the right hand side of
\erf07 gives us back $\hilde g$ when composed with $\id_X\Oti\R_M$,
one only has to use that $\hilde g\cir(\id_X\Oti\id_A\Oti\R_M)$ ${=}\,
\hilde g\cir(\id_X\Oti m\Oti\id_{\Dot M})$ and the unit property of $\eta$.
\\
The second isomorphism is shown analogously.
\hfill\checkmark

\vskip.3em

Thus indeed there are in particular isomorphisms
  \be  A\ota M\cong M \qquad{\rm and}\qquad N\ota A \cong N  \labl1a
  \end{equation}
as objects in \calc.

\vskip.3em

In general the object $M\Ota N$ will not be an \A-module, though
of course we have the

\lemma
When $M$ and $N$ are both \A-bimodules, then the tensor product $M\ota N$ 
has the structure of a {\rm(}left as well as right\/{\rm)} \A-module.

\vskip.86em
\rm

The tensor category of \A-bimodules plays a role in other circumstances 
(see e.g.\ \cite{brug1}), but it does not seem to be relevant in the context 
of boundary conditions. Rather, we like to restrict our attention even to 
a special class of \A-bimodules, namely those 
with a suitable compatibility condition between the left and right \rep. 
Such a class is e.g.\ naturally present for commutative algebras over the 
field of complex numbers. In that case every left module carries at the same
time the structure of a right module, and the two actions commute trivially. 
However, commutativity by itself is a too strong condition. Indeed, the same 
statement applies to supercommutative algebras, provided that the appropriate 
braiding in the category of super vector spaces (`rule of signs') is taken 
into account. In a sense, what we require is a generalization of the 
associativity \erf as of $m$, in which the middle factor of \A\ is replaced by 
an arbitrary \A-bimodule; in the \cft\ context, this should encode (part of)
the associativity of chiral operator product expansions.

Ideally we should deal with this issue by imposing additional restrictions
that can be formulated as properties of the \alg\ \A\ itself.
Here we use the following concept, formulated in terms of the category \calc,
that has been introduced in \cite{brug+}:\,%
 \footnote{~We thank Alain Brugui\`eres for providing us with
 information from \cite{brug+} prior to publication.} 

                    \vfill

\definition \labl sw
(i)\ For $A\iN\obj(\calc)$, an \A-{\sl \swap\/} (or briefly, {\sl \swap\/})
is a family of isomorphisms
$\,\wilde c_X\iN$\linebreak[0]$\Hom(X\Oti A,A\Oti X)$, one for each 
object $X\iN\obj(\calc)$, that obeys
  \be  (m\oti\id_X) \circ (\id_A\oti\wilde c_X) \circ (\wilde c_X\oti\id_A)
  = \wilde c_X \circ (\id_X\oti m)  \labl+b  \end{equation}
and
  \be  \wilde c_X \circ (\id_X\oti \eta) = \eta \oti \id_X  \labl+1
  \end{equation}
for all $X\iN\obj(\calc)$,
and that is tensorial and functorial, i.e.\ satisfies
  \be  \wilde c_{X\otimes Y}
  = (\wilde c_X \oti\id_Y) \circ (\id_X\oti\wilde c_Y)  \end{equation}
for all $X,Y\iN\obj(\calc)$ and
  \be  \wilde c_Y \circ (f\oti g) = (g\oti f)\circ\wilde c_X  \end{equation}
for all $f\iN\Hom(X,Y)$ and $g\iN\Hom(A,A)$.
\\[.3em]
(ii)\ Given an \A-\swap\ $\wilde c_X$, an \alg\ \A\ in \calc\ is called 
{\sl \swapc\/} iff it is $\tilde c$-com\-mu\-ta\-ti\-ve with
$\tilde c\eq\wilde c_A$\,.

                    \vfill

\remark
1.\ When \calc\ is braided, then the braiding $c_{A,X}$
and its inverse provide obvious candidates for the isomorphisms $\wilde c_X$. 
But even in a braided category other choices of $\tilde c$ may be of interest.
Moreover, typically the braiding and its inverse yield inequivalent \swap s,
so that even in the braided case deciding to
choose one of the two families of morphisms amounts to endowing \calc\ 
with additional structure. In particular, in general one cannot restrict to 
a subcategory for which the braiding and its inverse coincide without
loosing all information about symmetry breaking \bc s.
\\             [.3em]   
2.\ In the same way as the Frobenius property of and \alg\ \A\ corresponds to
s-t-duality in \cft\ (see remark \ref{st}.7), \swap-(co-)com\-mu\-ta\-ti\-vity
combined with the Frobenius property corresponds to what is known as
s-u-duality of four-point conformal blocks in the vacuum sector.
\\             [.3em]
3.\ When \A\ is \swapc\ and \swap-cocom\-mu\-ta\-ti\-ve, then the product 
\erf fg on $\Hom(A,A)$ is commutative.

                    \vfill 

\lemma
When \A\ is \swapc, then every left \A-module in \calc\ is naturally also 
a right \A-module, with the left and right actions of \A\ commuting.

\proof\dd
Let $N\eq(\Dot N,\R)$ be a left \A-module and set
  \be  \hat\R := \R\circ\wilde c_{\Dot N} \,.  \labl hr  \end{equation}
By using functoriality of the \swap, then the (left) \rep\ property of $\R$ 
and then the compatibility \erf+b between product and \swap, one sees that
  \be  \hat\R \circ (\hat\R \oti \id_A)
  = \hat\R \circ \llb \id_{\Dot N} \oti (m\cir\wilde c_A) \lrb
  \,.  \end{equation}
By \swap-commutativity this becomes $\hat\R\cir(\id_{\Dot N}\oti m)$,
and hence $\hat N\,{:=}(\Dot N,\hat\R)$ is a right \A-module.
Further, again by functoriality of the \swap\ and the \rep\ property
of $\R$ we have
  \be  \hat\R \circ (\R \oti \id_A)
  = \R \circ \llb (m\cir\wilde c_A) \oti \id_{\Dot N} \lrb
  \circ (\id_A \oti \wilde c_{\Dot N}) \,,  \end{equation}
which by using \swap-commutativity and applying the \rep\ property
of $\R$ backwards can be rewritten as $\R\cir(\id_A\oti\hat\R)$\,. 
\hfill\checkmark

\vskip.53em

{}From now on all left modules over a \swapc\ \alg\ \A\ will also be
regarded as right modules via \erf hr. Moreover, all right \A-modules
that we encounter are to be thought of as coming from left modules
via this construction.
Accordingly, it is preferable to reformulate the definition
of $M\Ota N$ using only left module structures. This is done as 
follows. First, as an object we have
  \be  \bearll &  M\ota N := \coker(f_1^{M,N}{-}\,f_2^{M,N}) \,,
  \\{}\\[-.5em]
  {\rm with}\quad & f_1^{M,N} := \hat\R_M \oti \id_{\Dot N}
  \equiv (\R_M\circ\wilde c_{\Dot M}) \oti \id_{\Dot N} \,,
  \quad f_2^{M,N} := \id_{\Dot M} \oti \R_N  \,, \eear  \labl12  \end{equation}
where both $\R_M$ and $\R_N$ refer to the left action. Then the \rep\ morphism
  \be  \R_{M\otimes_{\!A} N}^{}\;\in\Hom(A\Oti(M\Ota N),M\Ota N)
  \labl rm  \end{equation}
is the unique morphism that obeys $\R_{M\otimes_{\!A} N}^{}\cir
(\id_A\oti p_{M,N})\eq\hilde{\R}_{M\otimes_{\!A} N}$, with $p_{M,N}$ 
the projection \erf pc to the cokernel and (say) 
  \be  \bearl  \hilde{\R}_{M\otimes_{\!A}N}
  := p_{M,N}\circ f_1^{M,N}\cir ((\wilde c_{\Dot M}{)}^{-1}\Oti\id_{\Dot N})
  = p_{M,N}\circ(\R_M\oti\id_{\Dot N}) \\{}\\[-.7em] \hsp{17.7}
  \in \Hom(A\Oti\Dot M\Oti\Dot N,M\Ota N) \,.  \eear  \end{equation}
(Using that the left and right \rep\ of \A\ on $M$ commute, one checks that 
$\hilde{\R}_{M\otimes_{\!A} N}\cir(\id_A\oti(f_1^{M,N}\!{-}f_2^{M,N}))\eq0$.)
The left module properties of $M\Ota N$ then follow directly from
the left module structure of (say) $M$.

\vskip.3em

\lemma
The relations {\rm(}\ref{01}\,{\rm)} and {\rm(}\ref{02}\,{\rm)} are valid 
as isomorphisms of \A-modules.
As a consequence, one has the equalities, as \A-modules,
  \be  A\ota M = M \qquad{\rm and}\qquad N\ota A = N  \labl1c
  \end{equation}
for all $M,N\iN\obj(\calca)$ and
  \be  \F(X)\ota \F(Y) = \F(X\Oti Y)  \labl1e  \end{equation}
for all $X,Y\iN\obj(\calc)$.

\proof\dd
(i)~We only verify the statement for relation \erf02; the proof for \erf01
is analogous (recall that via \erf hr we can regard right modules as left 
modules). Thus we must show that the module structure of $N\Ota\F(X)$ that 
is defined by formula \erf rm coincides with the one provided by the
module structure of $N$ tensored with $\id_X$. Now by definition the 
former \rep\ is the unique morphism $\R_{N\otimes_{\!A}\F(X)}$ such that
$\R_{N\otimes_{\!A}\F(X)}\cir(\id_A\oti p_{N,\F(X)})\eq
\hilde\R_{N\otimes_{\!A}\F(X)}$ for
  \be  \hilde\R_{N\otimes_{\!A}\F(X)} = p_{N,\F(X)} \circ (\R_N \oti 
  \id_{A\otimes X}) = (\hat\R_N \oti \id_X) \circ (\R_N \oti \id_A \oti \id_X)
  \,.  \end{equation}
On the other hand, using that the two actions $\R_N$ and $\hat\R_N$ commute, 
one sees that the morphism $\R_N\Oti\id_X$ indeed does possess this property:
  \be  \bearll  (\R_N\oti\id_X) \circ (\id_A\oti p_{N,\F(X)}) \!\!\!\!
  &= \llb \R_N \circ (\id_A\oti\hat\R_N) \lrb \oti \id_X
  \\{}\\[-.8em]
  &= \llb \hat\R_N \circ (\R_N\oti\id_A) \lrb \oti \id_X
   = \hilde\R_{N\otimes_{\!A}\F(X)} \,.  \eear  \end{equation}
Thus the assertion follows by uniqueness.  
\\
(ii)~While the formulas \erf01 and \erf02 provide us a priori only with 
isomorphisms, by the arguments given in (i) we can invoke the fact that 
cokernels are unique only up to isomorphism to require that we even have 
equality
  \be  \F(X)\ota M = X\oti M  \qquad{\rm and}\qquad 
  N \ota \F(X) = N\oti X  \end{equation}
as \A-modules, with $\R_{X\otimes M}\eq\id_X\Oti\R_M$ and $\R_{N\otimes X}
\eq\R_N\Oti\id_X$. The relations \erf1c and \erf1e follow as special cases.
\hfill\checkmark

\vskip.42em

Based on the results above we can now show:

\proposition
When \A\ is \swapc, then \calca\ is a \tc\ with unit object \A.

                  \vfill
\proof\dd
(i)~According to formula \erf1c \A\ is a unit for the tensor product $\Ota$,
and the triangle identities are fulfilled trivially.
\\
(ii)~Making use of the defining properties of cokernels, by diagram
chasing -- the starting point being the observation that 
  \be  \bearl  \llb \id_M\Oti(f_1^{N,P}\!{-}f_2^{N,P}) \lrb \cir
  \llb (f_1^{M,N}\!{-}f_2^{M,N})\Oti\id_A \Oti\id_P \lrb  \\{}\\[-.5em] \hsp{9}
  = \llb (f_1^{M,N}\!{-}f_2^{M,N})\Oti\id_P \lrb \cir \llb \id_M\Oti\id_A\Oti
  (f_1^{N,P}\!{-}f_2^{N,P}) \lrb  \eear  \end{equation}
holds as a consequence of that fact that the left and right actions of \A\ on
$N$ commute -- for any triple $M,N,P\iN\obj(\calca)$ one constructs an 
isomorphism between $(M\Ota N)\Ota P$ and $M\Ota(N\Ota P)$. Verifying the
pentagon axiom is then another (lengthy) exercise in diagram chasing.
In the sequel we assume that (just like already in the case of the
underlying category \calc) the module category has
been traded for an equivalent strict tensor category by
invoking the coherence theorems, so that even
  \be  (M\ota N) \ota P = M \ota (N\ota P) \,.  \end{equation}
(iii)~Next we must define the tensor product $f\Ota g\iN\Hom(M\Ota N,
M'\Ota N')$ for every pair
of morphisms $f\iN\Hom_A(M,M')$ and $g\iN\Hom_A(N,N')$. It must be a
morphism such that $(f\Ota g)\cir p_{M,N}\eq\hilde{(f\Ota g)}$ for
  \be  \hilde{(f\Ota g)} := p_{M',N'}^{} \circ (f\Oti g)
  \; \in \Hom(M\Oti N,M'\Ota N') \,.  \end{equation}
Using the fact that $f$ and $g$ are morphisms of \calca\ together with the 
cokernel property of $p_{M',N'}$, one shows that $\hilde{(f\Ota g)}\cir
(f_1^{M,N}\!{-}f_2^{M,N})\eq0$. Together with the cokernel property of 
$p_{M,N}$ this in turn implies that a morphism $f\Ota g$ with the desired 
property exists and is unique.
\\
(iv)~Finally, one easily checks that the morphisms defined in (iii) satisfy
the relation $(f'\Ota g')\cir(f\Ota g)\eq(f'\cir f)\Ota(g'\cir g)$ whenever 
$f'\Oti g'$ and $f\Oti g$ are composable. We have thus verified all 
properties that a tensor category must possess.  
\hfill\checkmark

\remark
If \A\ is \hapl, then according to formula \erf1j the ground ring of \calca\
is just the ground ring $k$ of \calc, i.e.\ induction preserves the
ground ring.

Further, we have:

\proposition
When \A\ is \swapc, then the induction functor \F\ is a tensor functor.

\proof\dd
Tensoriality of $\F$ on objects just amounts to the identities $\F(\one)\eq A$
and \erf1e. To show tensoriality on morphisms we first note that, for all
$f\iN\Hom(X,X')$ and $g\iN\Hom(Y,Y')$ and all $X,X',Y,Y'\iN\obj(\calc)$,
the tensor product morphism $\F(f)\Ota\F(g)\iN\Hom_A(\F(X)\Ota\F(Y),\F(X')
\Ota\F(Y'))$ is by definition the unique morphism satisfying $(\F(f)\Ota\F(g))
\cir p_{\F(X)\otimes_{\!A}\F(Y)}\eq\hilde{(\F(f)\Ota\F(g))}$ with
  \be  \hilde{(\F(f)\Ota\F(g))} := p_{\F(X')\otimes_{\!A}\F(Y')} \circ 
  (\id_A\Oti f\Oti\id_A\Oti g) \,.  \end{equation}
On the other hand, using $p_{\F(X)\otimes_{\!A}\F(Y)}\eq\llb(m\Oti\id_X)\cir
(\wilde c_X\Oti\id_A)\lrb\Oti\id_Y$ and the analogous expression for $p_
{\F(X')\otimes_{\!A}\F(Y')}$, as well as functoriality of the \swap\ (applied
to $f$), one can check that the morphism $\id_A\Oti f\Oti g$ possesses this 
property. Thus $\F(f)\Ota\F(g)\eq\id_A\Oti f\Oti g\eq\F(f\Oti g)$.
\hfill\checkmark

\vskip.3em
On the other hand, the restriction functor
\G\ is not tensorial on objects; indeed, in general
$\G(M\Ota N)$ is only a quotient of $\Dot M\Oti\Dot N\eq\G(M)\Oti\G(N)$.  

\vskip.3em

\subsection{Dualities}
{}\mbox{$\ $}\\[.4em]
We now would like to establish left and right duality structures on
the module category \calca. We start with the

\definition  \labl51
When \calc\ possesses a (right) duality and an \A-\swap\ and \A\ is
\swapc, then for a left 
\A-module $M\iN\obj(\calca)$ the (right) {\sl dual module\/} $M\Vee$ is
  \be  M\Vee \equiv (\Dot M,\R)\Vee := (\Dot M\Vee\!{,}\,\R\Wee) \,,
  \end{equation}
where the morphism $\R\Wee$ is defined as
  \be  \bearl  \R\Wee \,{:=}\, (d_{\Dot M}\oti\id_{\Dot M\Vee}) \cir
  (\id_{\Dot M\Vee}\oti\R\oti\id_{\Dot M\Vee}) \cir
  ((\wilde c_{\Dot M\Vee}{)}^{-1}\oti b_{\Dot M}) \\{}\\[-.8em]\hsp{20}
  \,\in\Hom(A\Oti\Dot M\Vee\!{,}\,\Dot M\Vee) \,.  \eear  \end{equation}

In words, by making use of the presence of a \swap\
one is able to obtain the dual module by using the duality in \calc\
to deal with $\Dot M$ in place of $\Dot M\Vee$ and then employing the module 
structure of $M$. That $M\Vee$ is a left \A-module follows by combining
the \rep\ property of $\R$, the properties \erf+b and \erf+1 of
$(\wilde c_{\Dot M\Vee}{)}^{-1}$ and the duality axiom of \calc. The \swapy\ 
of \A\ is needed because use of the duality reverses the order in which the
two morphisms $\R$ are applied.

\vskip.42em

The following observation extends lemma \ref{0f}:
\lemma  \labl ff
A \tcc\ Frobenius \alg\ \A\ in a category with duality satisfies
$A\,{\stackrel\Phi\cong}\,\Av$ with \,$\Phi\iN\Hom_A(A,\Av)$. 
Thus \A\ is self-dual as an \A-module. 

\proof\dd
When \A\ is Frobenius, an isomorphism $\Phi\iN\Hom(A,\Av)$ is provided by 
formula \erf4f.
This is even a morphism of \calca, as is seen by using the functoriality
of $\wilde c$, associativity, and \tcy\ (as well as the explicit form of the 
\A-module structure of $\Av$ as given via definition \ref{51}). Likewise,
the inverse $\Phi^{-1}$ of $\Phi$ is an element of $\Hom_A(\Av\!{,}\,A)$,
so $\Phi$ is an {\em iso\/}morphism of \A-modules.
\hfill\checkmark

\vskip.63em

Analogously, when \calc\ also has a left duality, we have $\tilde\Phi\iN
\Hom_A(A,\Av)$, too, with $\tilde\Phi$ as defined by \erf ft. But by the 
\hapy\ \erf1i of \A\ the morphism set $\Hom_A(A,\Av)\,{\cong}\,\Hom(\one,\Av)$
is a free $k$-module of rank 1, and hence $\tilde\Phi$ differs from $\Phi$
only by a scalar.  As a consequence, the isomorphisms $\Nu_{\!A}^{(\tilde\Phi)}
\eq\tilde\Phi\cir\Phi^{-1}_{}$ \erf tp and $\tilde\Nu_{\!A}^{(\Phi)}
\eq\Phi\cir\tilde\Phi^{-1}_{}$ \erf pt satisfy $\Nu_{\!A}^{(\tilde\Phi)}
\eq\Nu_{\!A}^{(\Phi)}$ and $\tilde\Nu_{\!A}^{(\Phi)}\eq\tilde\Nu_{\!A}
^{(\tilde\Phi)}$ and are scalar multiples of $\id_{A\Vee}$, even though
\A\ is not \absi\ as an object in \calc.
These scalars can be calculated as in remark \ref{35}:
  \be  \Nu_{\!A}^{(\Phi)} = (\dimL(A))^{-1}_{}\bet\gam\, \id_{A\Vee} \,,
  \qquad \tilde\Nu_{\!A}^{(\Phi)} = (\dimR(A))^{-1}_{}\bet\gam\, \id_{A\Vee}
  \labl52  \end{equation}
with invertible $\diml(A)$ and $\dimr(A)$. Moreover, since the two morphisms 
are each other's inverses, we have $(\bet\gam)^2\eq(\dim(A))^2$, and hence 
  \be  \Nu_{\!A}^{(\Phi)} = \nu_{\!A}\,\id_{A\Vee}
  \qquad{\rm with}\qquad  \nu_{\!A} = \pm1 \,.  \labl na  \end{equation}
When combined with formula \erf a1, this in turn implies that
  \be  \dim(A) = \nu_{\!A}\,\bet\gam \,.  \labl a2  \end{equation} 

Further, together with identity \erf8x it follows that
  \be  d_A \circ (\id_{\!A} \oti (\eta\cir\eps)) \circ \tilde b_A 
  = \gam\,\nu_{\!A} \,.  \labl9x  \end{equation}
Similarly one shows, using also \tcy\ and the Frobenius axiom, that
  \be  \bearl  (d_A \oti \id_{\!A}) \circ (\id_{\!\Av} \oti \tilde c) \circ 
  (\tilde b_A \oti \id_{\!A}) \\{}\\[-.8em] \hsp6
  = ((\eps\cir m) \oti \id_{\!A}) \circ 
  (\id_{\!\Av} \oti \tilde c) \circ ((\Delta\cir\eta) \oti \id_{\!A})
  = \nu_{\!A}\, \id_{\!A} \,,  \labl9y  \eear  \end{equation}
as well as 
  \be  (\id_{\!A} \oti \tilde d_A) \circ (\tilde c \oti \id_{\!\Av}) \circ
  (\id_{\!A} \oti b_A)  = \nu_{\!A}\, \id_{\!A} \,.  \labl9z  \end{equation}

\vskip.3em

The analogue of lemma \ref{ff} for the case that
\calc\ has a left rather than a right duality reads $A\,{\cong}\,\Eev A$
with an isomorphism provided by $\tilde\Phi$. Here the left dual module
is $\Eev M\eq(\Eev{\Dot M}{,}\,\eew\R)$ with
  \be  \hsp{-.9}  \eew\R := (\id_{\Dot M\Vee}\Oti\tilde d_{\Dot M})
  \cir (\id_{\Dot M\Vee}\oti\R\oti\id_{\Dot M\Vee}) \cir
  ((\wilde c_{\Dot M\Vee}{)}^{-1}\Oti\id_{\Dot M}\Oti\id_{\Dot M\Vee}) \cir
  (\tilde b_{\Dot M}\Oti\id_{\Dot M\Vee}) \,.  \labl ld  \end{equation}

               \vfill 

While we already introduced the notion of a (right or left)
dual module as an object of \calca, we do not yet have a duality, 
in the sense of definition \ref{du}, in \calca. We proceed with the

\definition \labl BD
For $M\eq(\Dot M,\R_M)$ an \A-module in a sovereign category with an
\A-\swap, the morphism $B_M$ is defined as 
  \be  \bearll B_M \!\!\! &:= \bet^{-1}\,E(p^{}_{M,M\Vee} \cir b_{\Dot M})
  \\{}\\[-.8em]
  &\,= \bet^{-1}\,\R_{M\otimes_{\!A}M\Vee}
  \circ \llb \id_A \oti (p^{}_{M,M\Vee} \cir b_{\Dot M}) \lrb
  \;\in\Hom_A(A,M\Ota M\Vee) \,,  \eear  \end{equation}
and the morphism $D_M\iN\Hom_A(M\Vee\Ota M,A)$ is defined as the unique
morphism such that $D_M\cir p_{M\Vee\!,M}\eq \hilde D_M$ with
  \be  \hsp{-.6}\bearl  \hilde D_M := \tf{d_{\Dot M}}
  = (\id_A\oti d_{\Dot M}) \circ (\id_A\oti\R_M\Wee\oti\id_{\Dot M}) \circ
  ((\Delta\cir\eta)\oti\id_{\Dot M\Vee}\oti\id_{\Dot M})
  \\{}\\[-.8em] \hsp{22.7}
  \in\Hom(\Dot M\Vee\Oti\Dot M,A) \,.  \eear \end{equation}

               \vfill 

That these definitions are appropriate is shown by the following
              \newpage

\proposition \labl pr
When \A\ is a \swapc\ special Frobenius \alg, then the families of morphisms 
$B_M$ and $D_M$, with $M\iN\obj(\calca)$, furnish a duality on \calca.
\\
The bifunctor $\,\Ota{:}\ \calca{\times}\,\calca\,{\to}\,\calca$ is two-sided
exact, and the Grothendieck group of \calca\ inherits a natural ring structure.

\proof\dd
(i)~By use of the Frobenius property, the duality axiom and the unit property
of $\eta$, as well as the definition of the dual \rep\ $\R_M\Wee$, one
shows that $\hilde D_M\cir(f_1^{M\Vee\!,M}\!{-}f_2^{M\Vee\!,M})\eq0$, so that
$D_M$ indeed exists and is unique. 
\\[.16em]
(ii)~We write $B_M$ as $B_M\eq p_{M,M\Vee}{\circ}\,B_M^\flat$. Then 
$B_M^\flat\eq\bet^{-1}(\R_M\oti\id_{\Dot M\Vee})\cir(\id_A\oti b_{\Dot M}))$.
Also, one can rewrite $\hilde D_M$ as
$\hilde D_M\eq(\id_A\oti d_{\Dot M})\cir(\wilde c_{\Dot M\Vee}\oti\id_{\Dot M})
\cir(\id_{\Dot M\Vee}\oti\tf{\id_{\Dot M}})$. Using successively the duality
in \calc, \swapy, the \rep\ property of $\R_M$, the Frobenius axiom, the unit
property of $\eta$ and again \swapy, as well as functoriality of the \swap, one
then shows that
  \be  (\id_{\Dot M}\oti\hilde D_M) \circ (B_M^\flat\oti\id_{\Dot M})
  = \bet^{-1}\,(\wilde c_{\Dot M}{)}^{-1} \circ (\id_A\oti\R_M) \circ
  (\Delta\oti\id_{\Dot M}) \,.  \end{equation}
Further, after composing this morphism with $p_{M,A}\eq\hat\R_M$ we can use
again the \rep\ property of $\R_M$, as well as the special property of \A, 
to deduce that
  \be  p_{M,A} \circ (\id_{\Dot M}\oti\hilde D_M) \circ 
  (B_M^\flat\oti\id_{\Dot M}) = \R_M = p_{A,M} \,.  \end{equation}
With the help of the general relation \erf gh (specialized to $g\eq B_M\Ota
\id_M)$ and $h\eq\id_M\Ota D_M$) it therefore follows that
  \be  (\id_M\ota D_M) \circ (B_M\ota\id_M) = \id_M \,.  \end{equation}
The other duality axiom\, $(D_M\Ota\id_{M\Vee})\cir(\id_{M\Vee}\Ota B_M)
\eq\id_{M\Vee}$\, is proven analogously,
and the remaining statements are then immediate consequences.
\hfill\checkmark

\vskip.3em

\lemma \labl ll
On morphisms, the duality in \calca\ reduces to the one of \calc.

\proof\dd
In \calca, the morphism $f\Wee$ dual to $f\iN\Hom_A(M,N)$ is the
morphism in $\Hom_A(N\Vee,M\Vee)$ defined by
  \be  \hsp{-.6} f\Wee \,{:=}\; (D_N\oti\id_{M\Vee}) \circ (\id_{N\Vee}\oti f\oti
  \id_{M\Vee}) \circ (\id_{N\Vee}\oti B_M) \;\in \Hom_A(N\Vee\!{,}\,M\Vee) \,.
  \end{equation}
Inserting the definition of $B_M$ and $D_M$ and using the fact that $f$ is
a morphism in \calca\ as well as the Frobenius axiom and functoriality of the
\swap, one shows that the corresponding morphism $\hilde{(f\Wee)}$, satisfying
$p_{\!A,M\Vee}^{}{\circ}\,\hilde{(f\Wee)}\eq f\Wee$ is
  \be  \bearl  \hilde{(f\Wee)} = \bet^{-1} (d_{\Dot M} \oti \id_A \oti
  \R_{M\Vee}) \\{}\\[-.8em] \hsp{4.1}
  \circ \Llb \id_{M\Vee}  \circ \llb (f \oti \Delta \oti \id_{M\Vee}) \circ
  ((\wilde c_{\Dot M}{)}^{-1} \oti \id_{M\Vee}) \circ (\id_A \oti b_{\Dot M})
  \lrb \Lrb \,.  \eear \end{equation}
After composing with the projection $p_{A,M\Vee}$, the \rep\ morphism
$\R_{M\Vee}$ can be traded for a product $m$, and then the special property of 
\A\ cancels the prefactor $\bet^{-1}$. It follows that as a morphism in \calc,
i.e.\ in $\Hom(\Dot N\Vee\!{,}\,\Dot M\Vee)$, $f\Wee$ coincides with $f\Vee$ as
defined by using the duality in \calc.
\hfill\checkmark

\vskip.2em

Because of this result, in the sequel we will use the notation $f\Vee$ also
for the duality of morphisms in \calca. 

\vskip.3em

\proposition
When \A\ is Frobenius and \swapc, then the induction functor is compatible
with the duality, i.e.
  \be  (\F(X))\Vee \cong \F(X\Vee)  \labl vv  \end{equation}
for all $X\iN\obj(\calc)$ and
  \be  (\F(f))\Vee = \F(f\Vee)  \labl vw  \end{equation}
for all $X,Y\iN\obj(\calc)$ and $f\iN\Hom(X,Y)$.

\proof\dd
(i)~We have $\F(X\Vee)\eq(A\Oti X\Vee,\R_{A\otimes X\Vee})$ (with
$\R_{A\otimes X\Vee}\eq m\Oti\id_{X\Vee}$) and
$(\F(X))\Vee\eq((A\Oti X)\Vee,\R_{A\otimes X}\Wee)$.
An isomorphism $\phi_X$ between the objects $A\Oti X\Vee$ and $(A\Oti X)\Vee$
of \calc\ is provided by composing the \swap\ $(\wilde c_{X\Vee}{)}^{-1}{\in}
\,\Hom(A\Oti X\Vee\!{,}\,X\Vee\Oti A)$ with $\id_{X\Vee}\Oti\Phi\iN\Hom(X\Vee
\Oti A,X\Vee\Oti A\Vee)$ and the isomorphism between $X\Vee\Oti A\Vee$ and
$(A\Oti X)\Vee$. By lemma \ref{ff} and the definition of the
module structure of $\Av$ this is even an isomorphism in \calca.
Furthermore, using the duality in \calc, the Frobenius axiom, the properties
of $\eta$ and $\eps$ and \swapy, one checks that $\phi_X$ intertwines the 
action of \A, $\phi_X^{-1}\cir\R_{A\otimes X}\Wee\cir\phi_X\eq\R_{A\otimes X
\Vee}$. Thus $\F(X\Vee)$ and $(\F(X))\Vee$ are isomorphic as \A-modules.
\\[.1em]
(ii)~We have $(\F(f))\Vee\iN\Hom_A((\F(Y))\Vee,(\F(X))\Vee)$ and  
$\F(f\Vee)\iN\Hom_A(\F(Y\Vee)\!{,}\,\F(X\Vee))$, and we identify these morphism 
sets via the isomorphisms $\phi_X$ and $\phi_Y$ as constructed in (i). Then 
the equality \erf vw follows by similar manipulations as in the proof 
of lemma \ref{ll}.
\hfill\checkmark

\vskip.3em

Dealing with a left duality instead of a right duality of \calc\ works
in a completely analogous manner. In place of definition \ref{BD} we then
introduce duality morphisms $\tilde B_M$ and $\tilde D_M$ by
  \be  \bearll 
  \tilde B_M \!\!\! &:= \bet^{-1}\,E(p^{}_{M\Vee,\!M}\cir\tilde b_{\Dot M})
  \\{}\\[-.8em]
  &\,= \bet^{-1}\,\R_{M\Vee\otimes_{\!A}M}
  \circ \llb \id_A \oti (p^{}_{M\Vee\!,M} \cir b_{\Dot M}) \lrb
  \;\in\Hom_A(A,M\Vee\Ota M)  \eear  \end{equation}
and 
  \be  \hsp{-1.06}\bearl  \hilde{\tilde D_M} := \tf{\tilde d_{\Dot M}}
  = (\id_A\oti\tilde d_{\Dot M})\circ (\id_A\oti\R_M\oti\id_{\Dot M\Vee}) \circ
  ((\Delta\cir\eta)\oti\id_{\Dot M}\oti\id_{\Dot M\Vee})
  \\{}\\[-.8em] \hsp{23.4}
  \in\Hom(\Dot M\Oti\Dot M\Vee\!{,}\,A) \,.  \eear \end{equation} 
These indeed yield a left duality on \calca\ if \A\ is \swapc\ and
special Frobenius. Thus when \calc\ is autonomous, under the same conditions 
on \A\ the module category \calca\ is autonomous as well.

To derive the analogous statement also for sovereignty requires yet another 
condition, namely to select one of the two possible values \erf na
of the FS indicator:

\proposition
When \calc\ is sovereign and \A\ is a \swapc\ special Frobenius \alg\
with FS indicator $\nu_{\!A}\eq1$, 
then \calca\ is sovereign.

\proof\dd
(i)~By the definition of the left and right dual modules $\Eev M$ and $M\Vee$ 
(see formula \erf ld and definition \ref{51}) and
by sovereignty of \calc, they are identical as 
objects in \calc: $\G(\Eev M)\eq\Eev{\Dot M}\eq\Dot M\Vee\eq\G(M\Vee)$.
\\ 
(ii)~Using the duality axioms of \calc, the dual \rep s $\eew\R$ and 
$\R\Wee$ can be written in terms of the 
morphisms that are dual to $\R$ considered as a morphism in \calc:
  \be  \bearl
  \R\Wee = (\id_{\Dot M\Vee} \oti d_A) \circ (\R\Vee \Oti\, \id_A)
  \circ (\wilde c_{\Dot M\Vee}{)}^{-1} \,,
  \\{}\\[-.8em]
  \eew\R = (\id_{\Dot M\Vee} \oti \tilde d_A) \circ 
  (\wilde c_{\Dot M\Vee}{)}^{-1} \circ (\id_A \oti \Eev\R) \,.
  \eear \end{equation}
Inserting $\Eev\R\eq\R\Vee$, which holds by sovereignty of \calc, and
using functoriality of the \swap\ it follows that
  \be  \eew\R = \R\Wee \circ \llb ( (\id_{\!A} \oti \tilde d_A) \circ
  (\tilde c \oti \id_{\!\Av}) \circ (\id_{\!A} \oti b_A) ) \oti \id_{\Dot M\Vee}
  \lrb \,.  \end{equation}
By the identity \erf9z, the right hand side equals $\nu_{\!A}\R\Wee$, and
hence $\R\Wee$ iff $\nu_{\!A}$ is 1.
\\[.1em]
(iii)~Equality $\,\Eev f\eq f\Vee$ of left and right dual morphisms 
is an immediate consequence of lemma \ref{ll}.
\hfill\checkmark

\vskip.53em

One can then define left and right traces of morphisms $f\iN\Hom_A(M,M)$ by
  \be  \bearl
  \Trl(f) :=  D_M \circ (\id_{M\Vee}\ota f) \circ \tilde B_M \,, 
  \\{}\\[-.7em]
  \Trr(f) := \tilde D_M \circ (f\ota\id_{M\Vee}) \circ B_M \,,
  \labl Tr  \eear  \end{equation}
and left and right dimensions $\Diml(M)$ and $\Dimr(M)$ as the traces of 
$\id_M$. The results obtained above then imply the following

\lemma
When \calc\ is sovereign and \A\ is a \swapc\ special Frobenius \alg\
with $\nu_{\!A}\eq1$, then the traces defined in {\rm(}\ref{Tr}\,{\rm)} obey
  \be  \Trl(f) = \trl(f) \,/\, \dim(A) \,, \qquad
  \Trr(f) = \trr(f) \,/\, \dim(A) \,.  \end{equation} 
In particular,
  \be  \Diml(M) = \diml(M) \,/\, \dim(A) \,, \quad\,
  \Dimr(M) = \dimr(M) \,/\, \dim(A) \,.  \end{equation}

\proof\dd
Because of the identity \erf gh we may work just with morphisms in \calc.
By similar manipulations as in the proof of proposition \ref{pr},
one then arrives at 
  \be   \Trl(f) = \bet^{-1} (\id_A\oti d_{\Dot M}) \circ
  (\id_A\oti\R_{M\Vee}\oti f) \circ (\Delta\oti\tilde b_{\Dot M})
  \ \in\Hom_A(A,A) \,.  \end{equation}
for $f\iN\Hom_A(M,M)$. Since \A\ is an \absi\ module, $\Trl(f)$ must be a 
multiple of $\id_A$. (In the statement of the lemma we suppress the replacement
of $\trl(f)$ by $E(\eta\cir\trl(f))$ that makes $\trl(f)\iN k\eq\Hom(\one,\one)$
into an element of $\Hom_A(A,A)$, since it just reduces to tensoring by 
$\id_A$.) The proportionality constant can be determined by taking the trace 
and using the first of formulas \erf88 as well as \erf a2 together with the 
fact that $\nu_{\!A}\eq1$; the result is as claimed.\\
Right traces are treated analogously.
\hfill\checkmark

\vskip.3em

Note in particular that $\Dimr(A)\eq\Diml(A)\eq1$, as is indeed needed for
the tensor unit of a sovereign category. More generally,
  \be  \Dimr(\F(X)) = \dimr(X)  \qquad{\rm and}\qquad
  \Diml(\F(X)) = \diml(X)  \end{equation}
for all $X\iN\obj(\calc)$.

\vskip.3em

As an immediate consequence we arrive at the

\proposition
When \calc\ is spherical and \A\ is a \swapc\ special Frobenius \alg\
with FS indicator $\nu_{\!A}\eq1$, then \calca\ is spherical.

\rm 

\vskip.3em
\remark
1.\ The ultimate goal in \cft\ is to construct all correlation functions.
The ability to reach this goal is closely related to the issue of whether
the relevant \tc\ can supply invariants for a suitable class of links in a 
suitable class of three-manifolds (see \cite{fffs3}). Any modular category
over $\complex$ does this job, but generalizations are possible
\cite{kirI14,lyub10}. A mandatory property of the \tc\ seems to be
sovereignty. As seen above, this implies that the FS indicator $\nu_{\!A}$
must be equal to 1. This condition means that the self-conjugate object
$A$ is {\sl real\/} in the appropriate sense. In case \A\ is a direct sum
of simple subobjects each of which is a `simple current'
(see section 7), this reality condition should reduce to demanding that
the simple currents belong to the `effective center' \cite{gasc,krSc}.
\\
2.\ The dimensions provide us with distinguished numbers $\Diml(M)$ and
$\Dimr(M)$ for every module $M$, and with just a single $\Dim(M)$ when
\calca\ is spherical. In the particular case of \bc s preserving 
the full chiral \alg, this number is known, in \cft, 
as the `boundary entropy' of the \bc\ corresponding to $M$ \cite{aflu}. 
The dimension provides a natural generalization of this concept to
symmetry breaking \bc s; it shares the property of the boundary entropy
to behave multiplicatively under the tensor product of \calca.

\vskip.53em

\subsection{Semisimplicity}
{}\mbox{$\ $}\\[.4em]
We now discuss an issue that we did not touch upon when introducing
the category \calc, namely semisimplicity. Recall that an 
object $X'\iN\obj(\calc)$ is called a {\sl subobject\/} of $X\iN\obj(\calc)$
iff there exist morphisms $\imath\iN\Hom(X',X)$ and $p\iN\Hom(X,X')$
such that $p\cir\imath\eq\id_{X'}$; it is a {\sl proper\/} subobject iff
it is neither isomorphic to $X$ nor co-final. A {\sl simple\/} object is an 
object that does not possess a proper non-zero subobject. We call an \A-module 
$(M,\R)$ a {\sl simple module\/} iff it is simple as an object of \calca. A 
{\sl semisimple category\/} is a category for which every object is projective.
Since every short exact sequence $0\,{\to}\,M\,{\to}\,N\,{\to}\,P\,{\to}\,0$ 
with projective $P$ splits, every object of a semisimple category can in fact 
be written as a direct sum of simple objects.

\vskip.1em

There are several circumstances under which every \absi\ object is simple; 
then in particular the tensor unit, which is \absi\ by definition, is simple.
For instance, in a semisimple category, absolutely simple implies 
simple because the endomorphisms of an object $X$ that is a direct sum
$X\eq X_1{\oplus}X_2$ contain elements of the form $\xi_1\id_{X_1}{+}
\,\xi_2\id_{X_2}$ with $\xi_1\,{\ne}\,\xi_2$.
Without assuming semisimplicity, the implication still holds if the ground 
ring $k$ does not have zero divisors; for, if the projector $\imath\cir p
\iN\Hom(X,X)$ is a scalar multiple $\xi$ of $\id_X$, then
$\xi^2\id_X\eq\imath\cir p\cir\imath\cir p\eq\xi\id_X$, so by absence of
zero divisors $p$ and $\imath$ are invertible, hence $X$ and $X'$ isomorphic.
The implication is also valid if all morphism spaces are free 
$k$-modules, hence in particular when $k$ is a field. In contrast,
the converse implication is in general not valid even when 
the category is semisimple and $k$ is a field.\,%
 \footnote{~An example where the implication does not hold is provided
 by real irreducible representations of a real
 compact Lie group. In that case the endomorphism ring of a simple
 object is a commutative division algebra, and hence can be either the
 real numbers, the complex numbers, or the quaternions.} 

\remark \labl97
When the category \calc\ is semisimple, then the results obtained in the 
previous sections immediately lead to the following statements.
\\
1.\ If $M\eq(\Dot M,\R)\iN\obj(\calca)$ and $\Dot M$ is a simple object of
\calc, then by lemma \ref{42}.5 $M$ is a simple \A-module. In general
the converse does not hold.
\\
2.\ When the \alg\ \A\ in \calc\ is special, then by remark \ref{st}.2 we have
a direct sum decomposition
  \be  A\oti A \cong A \oplus X  \labl ab  \end{equation}
for some $X\iN\obj(\calc)$. 
When \A\ is Frobenius, then the coproduct is a morphism in \calca\ (see
formula \erf da); in fact, $m$ and $\Delta$ are then even morphisms
of \A-bimodules (for the natural bimodule structures of \A\ and $A\Oti A$).
Using this isomorphism, the product $m$ is nothing but
the projection to the first summand in \erf ab.
\\
3.\ When \A\ is \swapc, then as an object in \calc\ the tensor product 
$M\Ota N$ is a direct summand of $M\Oti N$, i.e.\
  \be  M\oti N = (M\ota N)\, \oplus Y  \end{equation}
with some $Y\iN\obj(\calc)$.
\\
4.\ By applying the functoriality property of an \A-\swap\ to the injections 
and projections that define subobjects, one sees that,
analogously as for a braiding, the \swap\ is already fully determined
by giving it for all simple objects of \calc.

\vskip.3em

Let us now study under which conditions semisimplicity is inherited by 
the module category \calca.

\lemma  \labl98
When \calc\ is semisimple and \A\ is Frobenius, then every induced \A-module 
$\F(X)$ is projective.

\proof\dd 
Consider any epimorphism $f\iN\Hom_A(M,N)$ of \A-modules.
Given a morphism $\hilde g\iN\Hom_A(\F(X),N)$, we write $\hilde{\Dot g}\,{:=}
\,E^{-1}(\hilde g)\iN\Hom(X,\Dot N)$, with $E$ the isomorphism \erf1k.
When \A\ is Frobenius, then according to lemma \ref{*3} $f$ is still an
epimorphism as a morphism of \calc. Since $X$ is projective (by semisimplicity
of \calc), this in turn implies that there is a $\Dot g\iN\Hom(X,\Dot M)$
such that $f\cir\Dot g\eq\hilde{\Dot g}$. We set $g\,{:=}\,E(\Dot g)
\iN\Hom_A(\F(X),M)$. Then using formula \erf fh we have
  \be  f \circ g = E(f\cir\Dot g) = E(\hilde{\Dot g}) = \hilde g \,.
  \end{equation}
Since this construction works for any epimorphism $f$ and any $\hilde g$,
every induced \A-module $\F(X)$ is projective. 
\hfill\checkmark

\vskip.3em

As a consequence we have the

\proposition\label{25}
When \calc\ is semisimple and the \calc-\alg\ \A\ is special Frobenius, 
then the module category \calca\ is semisimple.

\proof\dd
By the lemmata \ref{98} and \ref{99}, every \A-module $M$ is a submodule of
a projective module. Since any subobject of a projective object in a category 
is projective, $M$ is a projective module. 
\hfill\checkmark

\vskip.3em

\remark
1.\ Recall that semisimplicity of \calca\ implies that every absolutely
simple object is simple. It follows in particular that the
tensor unit \A\ of \calca\ is simple if \A\ is special Frobenius.
\\
2.\ If every object of the semisimple category \calc\ is a {\sl finite\/}
direct sum of simple objects, then in the direct sum decomposition
of a module in \calca\ only finitely many simple modules can be present, too.
In \cft, this situation is realized when one considers \bc s which all
preserve a {\sl rational\/} sub-vertex operator \alg\ of the bulk chiral \alg.

\section{Conformal \bc s from module categories}

Let us now return to the relationship between conformal \bc s and
the category theoretical concepts and results of the previous sections.  
The physical problem 
we want to address is the \class\ of modular invariants for a given
chiral conformal field theory and of \bc s associated to each such a modular 
invariant. A chiral CFT is characterized by its chiral algebra, a 
vertex operator algebra. According to our general philosphy, we
work with the tensor category of its representations rather than with
representation theoretic objects themselves.

\vskip.2em
\noindent
The picture we would like to put forward is the following: 
\\[.5em]
\hsp{1.1}\fbox{\parbox{32.9em}{%
  The physical modular invariant torus partition functions for a chiral
  conformal field theory with representation category \calc\ are in one-to-one
  correspondence with haploid special Frobenius algebras \A\ in \calc\ that 
  have trivial FS indicator. The partition function is of `extension type'
  if $A$ is commutative with respect to the braiding of \calc.}}
\\[.6em]
Based on this conjecture, we are led to a concrete picture of \bc s as well:
\\[.5em]
\hsp{1.1}\fbox{\parbox{32.9em}{%
  The \bc s for a \cft, with a given modular invariant partition function 
  and associated \calc-\alg\ \A,
  correspond to the objects in the category \calca\ of $A$-modules in \calc.}}
\\[.6em]
Furthermore, we expect the following dictionary between
physical concepts and mathematical notions (and also that the tensor category  
\calc\ that is relevant for the applications is sufficiently well-behaved 
such that those notions make sense):
\\[.2em]
{\bf 1)}~When the algebra object \A\ corresponds to a modular invariant of 
extension type, then the `boundary entropy' \cite{aflu} of a \bc\
is equal to the dimension of the corresponding object in \calca.
\\[.2em]
{\bf 2)}~In the extension case the annulus multiplicities 
are the structure constants of the Grothendieck ring of \calca, and the 
operator products of boundary fields are given by the $6j$-symbols of \calca.
\\[.2em]
{\bf 3)}~With the help of the category \calca\ the 
construction of correlation functions that was achieved  in \cite{fffs2,fffs3}
for the charge conjugation modular invariant and symmetry preserving \bc s
can be generalized to other modular invariants and to symmetry breaking \bc s.
\\[.2em] 
{\bf 4)}~In particular, the underlying modular invariant can be reobtained 
from \A\ by constructing the one-point correlation function
of the vacuum field on the torus.

\vskip.3em

The requirements we suggest to impose on an algebra object $A$ in order that 
it corresponds to a physical modular invariant are all motivated by properties 
that the vacuum field and conformal blocks on the sphere with vacuum insertions 
should possess in any consistent \cft. We already mentioned that the Frobenius 
property encodes s-t-duality, and \hapy\ the uniqueness of the vacuum state. 
The first part of the special property ensures that the tensor unit of \calc\ 
is, via the unit and the counit, in a distinguished manner a subobject of $A$; 
in other words, the original \voa\ (in particular, the Virasoro VOA) is in a
distinguished way a sub-VOA of the extended chiral algebra. \Swapy, if present, 
expresses the (graded) commutativity of \voa s, which in turn is in essence a 
locality property (only finite order poles are present in any operator product).  
Finally, in the extension case requiring the FS indicator to be 1 amounts to 
guaranteeing that the vacuum sector can have a definite conformal weight zero.

\vskip.3em

Various aspects of the picture suggested above are still under investigation. 
Here we content ourselves with presenting the following comments. Let us start 
with the relation between algebra objects and modular invariants. Most of the 
known modular invariants of conformal field theories can be described in terms 
of a group $\mathcal G$ of so-called \cite{scya} simple currents $\J$ -- that 
is, simple objects of a modular tensor category that have dimension one -- 
together with an element of $H^2(\mathcal G,\complex^\times)$ \cite{krSc}, 
called `discrete torsion'. The obvious candidate for the algebra object \A\ in
\calc\ is $A\eq\bigoplus_{\J\in \mathcal G}\!\J$. In general, this object can 
be endowed in several inequivalent ways with the structure of an associative 
algebra. In the next section we will see that for finite $\mathcal G$ the 
algebra structures are classified by $H^2(\mathcal G,\complex^\times)$, and 
therefore precisely correspond to the possible choices of discrete torsion.

Modular invariants that are not of simple current type are rare in
CFT; they are therefore called exceptional.  For many exceptional invariants 
the associated algebra object is known, too; an example will be discussed 
in section 8. The picture
we propose unifies simple current and exceptional modular invariants.

In the particular case of modular invariants of extension type (which includes
both simple current and exceptional invariants), the algebra object $A$ can be
directly read off, as an object of \calc, from the modular invariant. Namely,
such modular invariants amount to an extension of the chiral algebra by the 
primary fields that are paired with the vacuum, and accordingly $A$ is obtained 
by regarding the extended chiral algebra as an object in \calc. (This was the 
basis for the determination of the `extending sector' in \cite{fuSc14}.) 
In the non-extension case the relation between the algebra object and
the modular invariant (as well as the annulus partition functions) is more
involved; work on making the relationship explicit is in progress.
 
\vskip.3em

Our proposal provides a conceptual framework for the study of other problems
as well. For instance, the problem of deformations and moduli spaces of 
\cfts\ should be posed as a question about the deformation theory of 
algebras in \calc. Put this way, this question can be expected to make
sense even when \calc\ is not modular. Ideally, one might
wish to take for \calc\ a category as big as the one of all 
modules over the Virasoro algebra at fixed central charge. It is needless
to say that this is, technically, beyond reach.

A consistency check on our picture consists in the fact that the particular
class of algebras we specify is appropriate to ensure that various structural 
properties of \calc\ are inherited by the module category \calca. The 
following list collects the most important of these structures, together with
the properties of \A\ that are needed in order that they propagate from
\calc\ to \calca: 
\\[.63em]\hsp{1.2}
  \begin{tabular}{lcccc}
                &Frobenius&~~special~~& \swapc\ & haploid \&\ $\nu_{\!A}\eq1$
  \\[.5em]
  abelianness   &   +     &     +     &         &          \\[.2em]
  semisimplicity&   +     &     +     &         &          \\[.2em]
  tensoriality  &         &           &   +\    &          \\[.2em]
  duality       &   +     &     +     &   +\    &          \\[.2em]
  sovereignty   &   +     &     +     &   +\    & +$\quad$ \\[.2em]
  sphericity    &   +     &     +     &   +\    & +$\quad$ \\[.2em]
  \end{tabular}

\vskip1.1em

Let us conclude with a word of warning. One should not expect to be able 
to associate a `good' algebra object \A\ in \calc\ to every modular invariant 
partition function. This is much welcome, because it is known already for 
quite some time that there exist modular invariant partition functions which 
satisfy all the usual conditions but nevertheless are unphysical because they
cannot appear in any consistent CFT (for examples, see e.g.\ \cite{fusS}).
Our proposal encodes additional structural information. In our opinion it
therefore has a far better chance to provide us with the physically relevant
classification problem.

\section{The group case}

A large class of examples for our construction is supplied by the following
situation. Let \calc\ be modular over $k\eq\complex$ and
consider a \calc-\alg\ \A\ whose simple subobjects all have integral 
dimension and trivial twist. Then by  Deligne's reconstruction
theorem \cite{deli} (or, when \calc\ is a *-category, by the analogous
theorem in \cite{doro3}), the tensor subcategory \calco\ of \calc\ whose simple
objects are those isomorphic to the simple subobjects of \A\ is equivalent
to the \rep\ category of a finite group $G$, with the isomorphism classes of
simple objects of \calco\ being in one-to-one correspondence with those
of irreducible complex $G$-modules. We can then treat the objects of \calco\ 
as \findim\ vector spaces over $k\eq\complex$\,. Morphisms of \calco\ 
are then in particular linear maps and hence fully determined by their
action on the bases of those vector spaces.

For brevity let us further specialize to the case where $G$ is abelian.
Then the simple subobjects of \A\ have dimension 1, so that
$\dim(A)\eq|G|$, and the tensor product of \calc\ endows \A\ with the 
structure of the abelian group $G^*$, which is isomorphic to $G$. In the 
physics literature, objects of dimension 1 are known as {\sl simple 
currents\/} \cite{scya} and have been studied intensively; they also play 
an important role in the structure theory of fusion rings \cite{FRke}.

Associativity of $m$ then amounts to a cocycle condition. It follows that
the inequivalent \alg\ structures on \A\ are in one-to-one correspondence with
the elements of the cohomology group $H^2(G^*\!{,}\,k^\times)$.
The \alg\ structure is commutative iff it corresponds to the trivial class in
$H^2(G^*\!{,}\,k^\times)$. Let us choose a basis $\{\ej\}$ of \A, such that
$A\eq\bigoplus_{\J\in G^*}\!X_\J$ with $X_\J\eq k\ej$, $X_1\eq\one$ and\,%
 \footnote{~The objects $k\ej\Oti k\ek$ and $k\e{\J\K}$ are isomorphic.
 When the Frobenius\hy Schur indicator of all $X_\J$ equals $+1$, then 
 without loss of generality we may assume that these objects are actually equal.
 Otherwise we must be more careful; but those cases are not interesting
 (see below), and therefore we suppress this complication in the sequel.}
  \be  X_\J\oti X_\K = X_{\J\K} \,.  \labl jk  \end{equation}
Then a representative of the class of commutative \alg\ structures is given by
  \be  m(\ej\Oti\ek) := \e{\J\K} \,, \qquad  \eta(z) := z \e1 
  \labl71  \end{equation}
with $z\iN k$. By direct calculation one checks that together with
  \be  \Delta(\ej) := \sum_{\K\in G^*} \e{\J\K} \oti \e{\K\Vee} \,, \qquad
  \eps(\ej) := \delta_{\J,1} \e1  \end{equation}
this furnishes a cocommutative special Frobenius algebra with parameters
  \be  \bet = \dim(A) = |G| \,, \quad\ \gam = 1 \,.  \end{equation}
Up to a rescaling $\Delta\,{\mapsto}\,\beta^{-1}\Delta,\,\eps
\,{\mapsto}\,\beta\,\eps$ this is the unique special Frobenius structure
compatible with \erf71. Note that while $\bet$ and $\gam$ individually 
depend on the normalization of the coproduct, their product
$\bet\gam$, which equals $\dim(A)$, does not. 

We also remark that for abelian $G$ the tensor product of \calc\ is endowed 
with a $G$-grading that depends only on isomorphism classes, with the grading 
of a simple object $X$ given by
  \be  g_X\iN G:\quad g_X(\J):= S_{X,\J}/S_{X,\one} \,,  \labl74
  \end{equation}
where $S$ is the modular S-matrix \cite{TUra} of the modular category. (In 
\cft\ this quantity is called the `monodromy charge of the object $X$ with 
respect to the simple current $\J$'.) The object \A\ is homogeneous with 
trivial grading, $g_A\eq1$, so that \calca\ inherits a $G$-grading from \calc.
The full subcategory of \calc\ whose objects are 1-graded is premodular, and 
when applied to that subcategory, our construction is nothing but the 
modularisation procedure of \cite{brug2,muge6}, which in turn in physics 
terminology amounts to a `simple current 
extension' \cite{scya,fusS6} (since the associated modular invariant is then 
of pure extension type). In this situation, in particular the category 
$\calcao$ induced from the 1-graded subcategory of \calc\ is again modular, 
and in fact should be just the \rep\ category of an extended VOA \aext. The 
category $\calcao$ is a full tensor subcategory of the module category \calca;
its objects correspond to the `local sectors' in the physics literature. 

It is worth noting that most of simple current theory can be generalized to
cases where the twist $\theta$ is allowed to be non-trivial on some of the
simple subobjects of \A, provided that one makes the simplifying assumption
that one can still find a `section of objects' for the simple current group
$G^*$, i.e.\ a collection of objects $X_\J$, one for each $\J\iN G^*$,
such that \erf jk is still valid, i.e.\ such that the fusion product is
reproduced exactly (and not just up to isomorphism) by the tensor product.
(A property which does not generalize to the case with non-trivial twist is
that the grading \erf74, which is still present, is inherited by \calca.  
Incidentally, as soon as \A\ has irreducible subobjects 
with non-trivial twist, our construction is {\em not\/}\,%
 \footnote{~The corresponding statement in \cite{fuSc14} only
 applies to situations where \A\ has trivial twist.}
a categorical version of the $\alpha$-induction of \cite{lore,xu3,boev123}.) 
With this assumption we have in particular
  \be  \Hom(X_\J\Oti X_\K,X_\K\Oti X_\J) \cong k\,\id_{X_\J\otimes X_\K}
  \labl kj  \end{equation}
for all $\J,\K\iN G^*$. The compatibility 
  \be  \theta_{A\otimes A}\eq c_{A,A}\cir c_{A,A}\cir(\theta_A\Oti\theta_A)
  \labl aa  \end{equation}
between twist and braiding then implies that the braiding between two members
of the section has the form
  \be  c_{X_\J,X_\K} = \xi_{\J,\K}\,\id_{X_\J\otimes X_\K} \qquad{\rm with}
  \qquad \xi_{\J,\K}\,\xi_{\K,\J} = \theta_{X_\J\otimes X_\K} /
  \theta_{X_\J}\theta_{X_\K} \,.  \end{equation}
In particular, $\xi_{\J,\J}$ is determined by the twist up to a sign. Also, 
by tensoriality of the braiding one has $\xi_{\J\J',\K}\eq\xi_{\J,\K}{+}
\xi_{\J',\K}$ and $\xi_{\J,\K\K'}\eq\xi_{\J,\K}{+}\xi_{\J,\K'}$.

Now when the twist is non-trivial for some $\J\iN G^*$, then the \alg\
structure of \A\ is not commutative with respect to the braiding. To cure 
this we want to `change the braiding', i.e.\ to consider instead of the 
braiding $c$ a \swap\ $\wilde c$ such that \A\ is \tcc. Because of \erf kj 
$\tilde c\eq\wilde c_A$ must be of the form
  \be  \tilde c = \sum_{\J,\K\in G^*} \zeta_{\J,\K}\,c_{X_\J,X_\K} 
  \end{equation} 
with scalars $\zeta_{\J,\K}$, and because of \erf aa these are restricted by
  \be  \zeta_{\J,\K}\,\zeta_{\K,\J} = 1 \,.  \end{equation}

When combined with lemma 7.6 of \cite{kios}, this implies immediately that 
the object $\one{\oplus}X_\ell$ of the \mtc\ of the $\mathfrak{sl}(2)$ WZW 
model at level $\ell$ does not possess any \tcc\ \alg\ structure when $\ell$
is odd (and hence that there is no simple current modular invariant at odd 
level). We expect that this observation generalizes as follows. For any simple 
current we have $(\theta_{X_\J})^{2N_\J}\eq1$ with $N_\J$ the order of $\J$, 
i.e.\ the smallest positive integer such that $(X_\J)^{\otimes N_\J}\eq\one$;
a \swap\ such that the algebra $\bigoplus_{\J\in G^*}\!X_\J$ is \swapc\
should exist iff even $(\theta_{X_\J})^{N_\J}\eq1$ for all $\J\iN G^*$. 

We conclude this section with an instructive application
of the results of \cite{kirI14}. Suppose we are given a `holomorphic'
conformal field theory, i.e.\ a chiral CFT with a single primary field,
and a finite group $G$ (possibly non-abelian) acting on it by 
automorphisms of the chiral \alg. We are interested in those boundary
conditions which preserve the corresponding orbifold subalgebra.
It is shown in \cite{kirI14} that in this situation the tensor category 
$\calc$ of the orbifold subalgebra is equivalent to the category of \findim\
\rep s of the double D$(G)$, and that the algebra object is $\mathcal F(G)$,
the algebra of functions on $G$. According to theorem 3.2.\ of
\cite{kirI14}, the boundary category is equivalent to the category of
$G$-graded vector spaces, and the Grothendieck ring is the group ring
of $G$. (Notice that it is non-abelian for non-abelian groups $G$.)
In particular, irreducible boundary conditions are in one-to-one
correspondence to group elements $g\iN G$. On the level of boundary
states, this is also easy to understand: Denoting by $|B\rangle
\iN\mathcal H\oti\mathcal H$ the unique Cardy boundary state on the single
irreducible module $\mathcal H$ of the holomorphic CFT, the symmetry breaking
boundary states are just $|B,g\rangle\,{:=}\,(R(g)\Oti{\bf1})|B\rangle$.  

\section{The $E_6$-type invariant of $A_1^{\scriptscriptstyle(1)}$}

As an example that goes beyond the group case of section 7, take \calc\ to 
be the \mtc\ associated to the $A_1^{\scriptscriptstyle(1)}$ Lie algebra
at level $\ell\eq10$, with 11 isomorphism classes $\lambda\iN\{0,1,2,
...\,,10\}$ of simple objects, and \A\ to be in the class $0{\oplus}6$.
This corresponds to the $E_6$-type modular invariant of
$A_1^{\scriptscriptstyle(1)}$ and has been studied from various
points of view \cite{prss,xu3,boev123,bppz2,fuSc14}.

The object $A\eq\one\,{\oplus}\,X_6$ is a commutative algebra, as follows 
\cite{kios} by the existence
of a conformal embedding $\abar\;{\hat=}\,{(A_1^{\scriptscriptstyle(1)})}
_{10}\,{\subset}\,{(B_2^{\scriptscriptstyle(1)})}_1\,{\hat=}\;\aext$.
Since \A\ is the direct sum of only two simple subobjects, it is also not
difficult to extend this structure to the one of a special Frobenius algebra;
we have not analyzed whether this can be done in a unique way.

The category \calca\ is known to possess six isomorphism classes of simple 
objects, which we label as $\oa,\,\va,\,\sa,\,\ob,\,\vb,\,\sb$. Via the 
forgetful functor $\G$ they correspond to isomorphism classes of \calc-objects as
  \be  \bearll
  \G(\oa) = [\A] = 0 \oplus 6 \,,\ & \G(\ob) = 1\oplus 5\oplus 7 \,,
  \\{}\\[-.6em]
  \G(\va) = 4 \oplus 10 \,,        & \G(\vb) = 3\oplus 5\oplus 9 \,,
  \\{}\\[-.6em]
  \G(\sa) = 3\oplus 7 \,,  & \G(\sb) = 2\oplus 4\oplus 6\oplus 8
  \,. \eear \labl95  \end{equation}
The fusion rules are commutative, with $\oa$ acting as unit element
and \cite{xu3,fuSc14}
  \be  \bearlll
  \va \star \va = \oa \,,   &  \va \star \ob = \vb \,, 
                            &  \ob \star \ob = \oa \oplus \sb \,,
  \\{}\\[-.76em]
  \va \star \sa = \sa \,,   &  \va \star \vb = \ob \,,
                            &  \ob \star \vb = \va \oplus \sb \,,
  \\{}\\[-.76em]
  \sa \star \sa = \oa\oplus\va \,, \ & \va \star \sb = \sb \,, 
                            &  \ob \star \sb = \sa \oplus \ob \oplus \vb \,,
  \\{}\\[-.76em]
  & \sa \star \ob = \sb \,, &  \vb \star \vb = \oa \oplus \sb \,, 
  \\{}\\[-.76em]
  & \sa \star \vb = \sb \,, &  \vb \star \sb = \sa \oplus \ob \oplus \vb \,, 
  \\{}\\[-.76em]
  & \sa\star\sb = \ob\oplus\vb \,,\ \ &  \sb\star\sb = \oa\oplus\va\oplus 2\sb
  \,. \eear  \end{equation}

In view of the decompositions \erf95 it is of interest to study the
transformation of the functions
  \be  \bearll
  \Chi_\oa := \chii_0+\chii_6  \,,     & \Chi_\ob := \chii_1+\chii_5+\chii_7 \,,
  \\{}\\[-.8em]
  \Chi_\va := \chii_4+\chii_{10}\,,\ \ & \Chi_\vb := \chii_3+\chii_5+\chii_9 \,,
  \\{}\\[-.8em]
  \Chi_\sa := \chii_3+\chii_7 \,,      &  
  \Chi_\sb := \chii_2+\chii_4+\chii_6+\chii_8  \eear  \end{equation} 
with $\chii_\lambda$ denoting the character of the irreducible highest weight
$A_1^{\scriptscriptstyle(1)}$-module of isotype $\lambda$,
under the modular group action $S{:}\;\tau\,{\mapsto}\,{-}1/\tau$.
Expectedly, while the functions $\Chi_a$ (which are in fact irreducible
characters for $B_2^{\scriptscriptstyle(1)}$ at level 1) transform into
each other, the space of all six functions is not closed under this
transformation. But, remarkably enough, one does get closure if one includes 
in addition the functions
  \be  \breve\Chi_\oa := \cald^{1/2}\chii_0 - \cald^{-1/2}\chii_6 \,, \quad\
  \breve\Chi_\va := \cald^{-1/2}\chii_4 - \cald^{1/2}\chii_{10}   \,, \quad\
  \breve\Chi_\sa := \chii_3 - \chii_7  \end{equation}
with
  \be  \cald := \dim(6) = 2+\sqrt3 \,. \quad \end{equation}
Indeed, we find
  \be  \Chi_a(\mbox{-$\frac1\tau$}) = \sum_b S_{a,b}\, \Chi_b(\tau)\,, \quad\;
  \Chi_{\check a}(\mbox{-$\frac1\tau$}) = \sum_b S_{a,b}\, \breve\Chi_b(\tau)
  \,, \quad\;
  \breve\Chi_a(\mbox{-$\frac1\tau$}) = \sum_b S_{a,b}\, \Chi_{\check b}(\tau)
  \end{equation}
with $a,b\iN\{\oa,\va,\sa\}$ and a single $3{\times}3$-matrix $S$ with entries
  \be  \Llb S_{a,b} \Lrb := \mbox{\large$\frac12$}\,\left( \begin{array}{ccc}
  1 & \sqrt2 & 1 \\[.12em] \sqrt2 & 0 & -\sqrt2 \\[.12em] \;1 & -\sqrt2 & \;1 
  \eear \right) ,  \end{equation}
which is nothing but the modular S-matrix of the extended theory $B_2$
at level 1. This structure is intriguingly close to the structure
of $\zet_2$-orbifold theories based on inner automorphisms (compare
e.g.\ \cite{bifs}). One is therefore tempted to look for an
analogue of the twining characters of the twisted sectors, which
transform into each other under the S-transformation. Indeed, the whole 
11-dimensional space of level 10 $A_1^{\scriptscriptstyle(1)}$-characters
is spanned by the nine functions $\Chi_a$, $\Chi_{\check a}$ and $\breve\Chi_a$
together with $\chii_5$ and
  \be  \breve\Chi_- := (\chii_0-\chii_2+\chii_4-\chii_6+\chii_8-\chii_{10})
  /\sqrt6  \end{equation}
which transform into each other as
  \be  \chii_5(\mbox{-$\frac1\tau$}) = \breve\Chi_-(\tau) \,, \qquad
  \breve\Chi_-(\mbox{-$\frac1\tau$}) = \chii_5(\tau) \,.  \end{equation}
The close similarity of this structure with orbifolds, \resp\ simple
current extensions, gives further evidence to the idea that exceptional
modular invariants can be understood through a suitable generalization of
the concepts encountered in orbifold theories.

\section{NIM-reps}
         
\definition
A {\sl NIM-rep\/} of the Grothendieck ring of a sovereign tensor category
is a representation $R$ by \findim\ matrices with non-negative integral
entries satisfying
  \be  R(X\Vee) = (R(X))^{\rm t}_{} \qquad{\rm and}\qquad
  R(\one) = \onematrix \,.  \labl1'  \end{equation}

This concept of a NIM-rep has appeared in the search for modular invariants
\cite{dizu}, in the analysis of annulus amplitudes in \cft\ (see e.g.\
\cite{sasT2,fuSc112,bppz2}), and in the study of the general structure of
fusion rings \cite{gann16}. It is worth noting, however, that the \class\ 
of NIM-reps is not the same as the \class\ of good algebra objects. Rather, 
a NIM-rep can be unphysical, i.e.\ not correspond to any consistent set 
of \bc s, similarly as there exist unphysical modular invariant \parfu s.

NIM-reps also naturally arise in the context of $\alpha$-induction
\cite{boev123}, in a way which can be easily formulated in category
theoretic terms. Since NIM-reps are representations by \findim\ matrices,
their natural habitat are rational \cfts. In our context this means that 
\calc\ is semisimple and has only finitely many isomorphism classes of simple
objects. By the results of subsection 5.4 these properties are then also
shared by the module category \calca, provided that \A\ is a special Frobenius 
\alg. In addition, we assume that \calca\ is a tensor category, which requires 
that \A\ is also \swapc.

\proposition
When $k$ is a field and \calc\ is sovereign, then the dimension matrices
  \be  (R(X))_m^n := \dim\,\Hom_A(\F(X)\Ota M,N)  \labl2'  \end{equation}
and
  \be  (\tilde R(X))_m^n := \dim\,\Hom_A(M\Ota \F(X),N) \,,  \labl2"
  \end{equation}
with rows and columns labelled by the isomorphism classes $m\eq[M]$
of simple \A-modules in \calca, furnish a right and a left NIM-rep of 
the Grothendieck ring of \calc, \resp.
\proof\dd
The fact that \calca\ is semisimple implies (see e.g.\ \cite{TUra}) the 
following {\sl domination property\/}: For every pair of modules $M,N\iN\obj
(\calca)$, every morphism $f\iN\Hom_A(M,M')$ can be written as a finite sum
  \be  f = \sum_{\ell=1}^{\ell_{\rm max}(f)} j_\ell \circ p_\ell  \end{equation}
with $p_\ell\iN\Hom_A(M,N_\ell)$ and $j_\ell\iN\Hom_A(N_\ell,M')$ for
suitable (not necessarily distinct) simple objects $N_\ell\iN\obj(\calca)$.
In particular,
  \be  \dim\,\Hom_A(M,M') = \sum_{[N]:\; N\;{\rm simple}}
  \dim\,\Hom_A(M,N) \cdot \dim\,\Hom_A(N,M') \,,  \end{equation}
where the sum extends over the equivalence classes of simple objects of
\calca. For brevity we will write such summations from now on just as
$\sum_{[N]}$, and analogously as $\sum_{[X]}$ when dealing with \calc.
It follows that the dimension matrices \erf2' satisfy
  \be \bearll   R(X\Oti Y)_m^{m'} \!\!\!
  &= \dim\,\Hom_A(\F(X\Oti Y)\Ota M,M')
  \\{}\\[-.6em] 
  &= \dim\,\Hom_A(\F(X)\Ota\F(Y)\Ota M,M')
  \\{}\\[-.46em]
  &= \dim\,\Hom_A(\F(Y)\Ota M,\F(X)\Vee\Ota M')
  \\{}\\[-.6em]
  &= \dsty\sum_{[N]} \dim\,\Hom_A(\F(Y)\Ota M,N) \; \dim\,\Hom_A(N,\F(X)\Vee\Ota M')
  \\{}\\[-.86em]
  &= \dsty\sum_{[N]} \dim\,\Hom_A(\F(Y)\Ota M,N) \; \dim\,\Hom_A(\F(X)\Ota N,M')
  \\{}\\[-.86em]
  &= \dsty\sum_{n} R(Y)_m^n \, R(X)_n^{m'} \,,
  \eear \end{equation}
where we used tensoriality of $\F$, duality in \calca,
domination in \calca, and again duality in \calca.
This is nothing but the defining property of a right \rep. 
\\
The left \rep\ property of $\tilde R$ follows analogously:
  \be \bearll   \tilde R(X\Oti Y)_m^{m'} \!\!\!
  &= \dim\,\Hom_A(M\Ota\F(X)\Ota\F(Y),M')
  \\{}\\[-.52em]
  &= \dim\,\Hom_A(M\Ota\F(X),M'\Ota\F(Y)\Vee)
  \\{}\\[-.6em]
  &= \dsty\sum_{[N]}\dim\,\Hom_A(M\Ota\F(X),N) \;\dim\,\Hom_A(N,M'\Ota\F(Y)\Vee)
  \\{}\\[-.86em]
  &= \dsty\sum_{n} \tilde R(X)_m^n \, \tilde R(Y)_n^{m'} \,.
  \eear \end{equation}
To show also the property \erf1', we can use the identity \erf vv 
together with $\dim\,\Hom_A(M,M')\eq\dim\,\Hom_A(M',M)$ to conclude that
  \be  \hsp{-.6} \bearll (R(X\Vee))_m^n   \!\!\!
  &= \dim\,\Hom_A(\F(X\Vee)\Ota M,N) = \dim\,\Hom_A((\F(X)\Vee\Ota M,N)
  \\{}\\[-.8em]
  &= \dim\,\Hom_A(M,\F(X)\Ota N) = \dim\,\Hom_A(\F(X)\Ota N,M)
  \\{}\\[-.8em]
  &= (R(X))_n^m \,,  \eear  \end{equation} 
and analogously for $\tilde R$. The second part of \erf1' is immediate:
  \be  (R(\one))_m^n = (\tilde R(\one))_m^n = \dim\,\Hom_A(M,N)
  = \delta_{m,n}^{} \,.  \end{equation}
\vsp\hfill\checkmark

\vskip.3em

Also note that 
  \be  \hsp{-.35} R(X)_m^n = \dim\,\Hom_A(\F(X)\Ota M,N) = \dim\,
  \Hom_A(N\Vee\Ota\F(X)\!{,}\,M\Vee) = \tilde R_{n^\vee}^{m^\vee}
  .  \end{equation}
Since the relabelling $M\,{\mapsto}\,M\Vee$ just constitutes a basis change
of the Grothendieck ring of \calca, this means that $R$ and $\tilde R$ are,
as abstract representations, contragredient to each other.
When the Grothendieck ring of \calca\ is commutative, the two \rep s are,
of course, identical.

Let us also explore what we can learn, conversely,  about
the \alg\ \A\ by knowing a NIM-rep $R$ that is furnished by the
expression \erf2' (or \erf2"). When applied to the diagonal matrix element
$R(X)_{[A]}^{[A]}$ that corresponds to the unit object $A$ of 
\calca, the reciprocity relation \erf4g yields
  \be  R(X)_{[A]}^{[A]} = \dim\,\Hom_A(\F(X),A) = \dim\,\Hom(X,\G(A)) \,.
  \end{equation}
Thus, a NIM-rep of the Grothendieck ring of \calc\ determines an underlying 
algebra \A\ (if it exists) uniquely, up to isomorphism, as an object in \calc:
  \be  A \,\cong\, \bigoplus_{[X]}  R(X)_{[A]}^{[A]} \, X \,.  \end{equation} 
To apply this formula we must, however, still know \A\ (\resp\ its isomorphism
class) as an object of \calca, or in other words, which element $[M_0]$ of 
the distinguished basis of the Grothendieck ring of \calca\ corresponds
to the tensor unit of \calca. In practice this is often far from clear;    
the following observation is therefore useful: For every module $N$ we have
  \be  R(X)_n^n = \dim\,\Hom_A(\F(X)\Ota N,N)
  = \dim\,\Hom_A(\F(X),N\Ota N\Vee) \,.  \end{equation} 
Since the tensor product module $N\Ota N\Vee$ contains the tensor unit $A$,
by domination we have $R(X)_n^n\,{\ge}\,R(X)_{[A]}^{[A]}$. As a consequence, 
from the NIM-rep the algebra can be determined as an object of \calc, even 
without knowing the basis element $[M_0]$ that corresponds to the algebra, as
  \be  A \,\cong\, \bigoplus_{[X]} \min_m R(X)_m^m\; X \,.  \labl mi
  \end{equation} 
Also, since $R(\one)$ is the unit matrix we have
$\dim\,\Hom(\one,A)\eq1$, and the identity \erf1' tells us that $\dim\,
\Hom(X,A)\eq\dim\,\Hom(X\Vee\!{,}\,A)$. Thus \A\ is \hapl\ and self-dual.

Now as already mentioned, the physically meaningful task is to classify 
algebra objects in \calc, but not every NIM-rep gives rise to a consistent 
algebra object. Indeed, in the situation considered here,
an associated algebra object can exist only if there is an $m_0$ such that
the minimum in formula \erf mi is attained at $m\eq m_0$ simultaneously
for all basis elements $[X]$ of the Grothendieck ring of \calc; when such
an $m_0$ does not exist, then the NIM-rep is necessarily unphysical.  
(Still, classifying NIM-reps can be a useful intermediate step; compare 
the strategy reviewed in \cite{zube8} to obtain an overview over possible
\bc s, or the use of the results of \cite{etkh}
in the classification \cite{kios} of algebra objects in representation
categories for $A_1^{\scriptscriptstyle(1)}$ at positive integral level.)
Note that in order to reach this conclusion we had to require \calc\ to 
possess a \swap. As we cannot say much about \swap s that do not coincide 
with the braiding, our result is of direct help only
when \A\ is commutative with respect to the braiding
and thus when the torus partition function is of extension type.

\section{Solitonic sectors and \voa s}

Let us recapitulate the situation we are studying: From the algebraic and 
algebro-geometric data of a chiral \cft\ -- a \voa, its category of 
representations and the associated system of conformal blocks (see e.g.\ 
\cite{BEfr,Scfw}) -- we have abstracted a tensor category \calc, and investigated 
the categories \calca\ of modules of suitable algebra objects in \calc.

An obvious question at this point is how to interpret the objects
of \calca\ at the algebraic level. Mathematically, this issue should
be compared to the following situation. The so-called
center-construction (see e.g.\ \cite{KAss}) associates to every
strict tensor category $\calc$ a braided tensor category $\mathcal Z(\calc)$.
This construction can in particular be applied to the tensor category
$\calc(H)$ of modules over a \findim\ Hopf algebra $H$ with invertible 
antipode.  On the other hand, at the algebraic level one can 
construct the quantum double D$(H)$ of $H$, whose representation category 
$\calc({\rm D}(H))$ is braided. The category theoretic 
and the algebraic construction commute in the sense that the categories
$\calc({\rm D}(H))$ and $\mathcal Z(\calc(H))$ are equivalent categories.

In the present context, the algebraic side of the construction remains open
for the moment. Concretely, the question is: Can the objects of a module 
category $\calca\eq {(\calcm)}_{\!A}$ be interpreted as generalized modules 
over some algebraic structure that extends the original vertex algebra \abar? 
This question is also imposed to us by
physical considerations: At the category theoretic level, the transition from
chiral \cft\ (living on a complex curve) to full \twodim\ \cft\ (living on
a real \twodim\ manifold which may be non-orientable and may have a
boundary) can be analysed model-independently, including the classification
of \bc s and computation of the model-independent part of correlation functions.
But still part of the physically interesting information is missing, in 
particular about the structure of the space of states.

Obtaining this information  
would be simpler if we had an answer to the following challenging
problems in the theory of vertex algebras: First, to characterize vertex 
algebras as algebra objects in suitable categories. (For work in this 
direction see \cite{borC22,soib4}). Second, to formulate a framework that
allows to reconstruct a vertex algebra from its representation theory, a fiber 
functor and possibly some additional data. Both constructions would further 
corroborate the status of vertex algebras as natural algebraic objects. 

\vskip.3em

So far we are not able to provide concrete prescriptions for the algebraic
side of the construction.\,%
 \footnote{~For simple current extensions of the \mtc\ that is associated to
 an affine Lie \alg\ at positive integral level \A\ is known \cite{dolm} to 
 have the structure of a so-called abelian intertwining \alg. But except when
 this happens to be even a \voa, it is not the structure we are looking for.} 
But the following independent observations shed some light on the situation:

First, at least in many concrete models the categories in question can be
realized in the setting of C$^*$- and von Neumann \alg s 
\cite{EVka,boev123,muge7}, where the objects of the boundary category 
correspond to endomorphisms of such \alg s. When those endomorphisms are 
localizable in a bounded domain, they give rise to distinguished \rep s which,
in turn, are the counterparts of the ordinary \rep s in the \voa\ 
setting. But there are also objects whose associated endomorphisms are 
{\sl solitonic\/}, i.e.\ \cite{froh3,fred10} only localizable in certain 
unbounded domains. These still correspond to \rep s of the relevant C$^*$-\alg,
but do not possess analogues among the ordinary \rep s of a \voa.

Second, when the sub\alg\ \abar\ of the \voa\ \aext\
that is respected by all \bc s corresponding to \calca\
is an orbifold sub\alg, i.e.\ is the fixed point set under the action of a
(finite) group $G$ of automorphisms of \aext, the solitonic \rep s are already 
known.  They are provided by twisted \rep s of \aext \cite{dolm3} and provide 
the twisted sectors of the $G$-orbifold of the extended \cft. In the category 
theoretic framework, this situation has recently been studied in \cite{kirI14} 
(the untwisted sector had been analyzed earlier in \cite{brug2,muge6}). It is 
shown in \cite{kirI14} that the notion of a twisted \rep\ of a \voa\ has a 
precise correspondence in the form of a {\sl twisted module\/} over a \calc-\alg\
\A. The definition of twisted modules, in turn, can be taken over rather 
directly to the situation where the twists do no longer form a group.

Third, for some classes of models certain vertex operator (super)algebras 
associated to simple currents have been constructed in \cite{lI5,femi,lI10}. 
In particular, in the super-case the objects correspond to 
super-representations. A feature that this observation has in common with the 
previous point is that
the \rep s are no longer defined in the category of vector spaces, but rather 
on super vector spaces or, in the case of twisted modules, $G$-vector spaces.

A fourth observation concerns the spaces of conformal blocks. It was already 
observed some time ago \cite{blma} that for simple currents of half-integral
conformal weight the factorization rules become much more transparent
once one promotes the vector spaces of conformal blocks to super vector spaces.
Thinking about conformal blocks as (generalized) co-invariants \cite{BEfr},
this nicely fits with our third remark. On the other hand, in the category 
theoretic approach the spaces of conformal blocks are constructed from
morphism spaces. The observation just mentioned therefore suggests to look
for suitable extensions of those spaces, compare remark \ref{42}.8.

Finally, it is apparent from the examples that we 
discussed above (and from many more) that the character sums
  \be  \chii_M := \hsp{-.3}\sum_{[X]:\;X\in\obj(\mathcal C)\;{\rm simple}}
  \hsp{-.8}\dim\,\Hom(X,\G(M))\,\chii_X  \end{equation}
possess interesting transformation properties under modular group, \resp\ 
under some subgroup of it. In the super-case, the relevant subgroup is 
precisely the one that respects spin structures on the torus.

\vskip.2em

Together these observations suggest that the objects of \calca\ should 
find their representation theoretic interpretation in a category that 
extends vector spaces, in a similar way as super vector spaces do.
A corresponding extension should exist for complex curves as 
well, very much like supersymmetric curves constitute an extension 
of complex curves, and such that a correct definition of conformal 
blocks as generalized vector spaces arises as a natural consequence. 

\vskip1.5em {\bf Acknowledgment}.
It is a pleasure to thank A.\ Brugui\`eres, T.\ Gannon, A.A.\ Kirillov Jr., 
G.\ Masbaum, M.\ M\"uger, I.\ Runkel, A.N.\ Schellekens and N.\ Sousa 
for helpful discussions and correspondence.  

 \newcommand\wb{\,\linebreak[0]} \def\wB {$\,$\wb}
 \newcommand\Bi[1]    {\bibitem{#1}}
 \renewcommand\J[5]   {{\em #5}, {#1} {\bf #2} ({#3}), {#4} }
 \newcommand\JL[4]    {{\bf #2} ({#3}), {#4}}
 \newcommand\JM[4]    {{#1} {\bf #2} ({#3}), {#4}}
 \newcommand\JN[4]    {{#1} {\bf #2} ({#3})}
 \newcommand\Prep[2]  {{\em #2}, pre\-print {#1}}
 \newcommand\PRep[2]  {{\em #2}, {#1}}
 \newcommand\BOOK[4]  {{\em #1\/} ({#2}, {#3} {#4})}
 \newcommand\inBO[7]  {{\em #7}, in:\ {\em #1}, {#2}\ ({#3}, {#4} {#5}), p.\ {#6}}
 \def\jf    {J.\ Fuchs}
 \def\adma  {Adv.\wb Math.}
 \def\aspm  {Adv.\wb Stu\-dies\wB in\wB Pure\wB Math.}
 \def\atmp  {Adv.\wb Theor.\wb Math.\wb Phys.}
 \def\coia  {Com\-mun.\wB in\wB Algebra}
 \def\coma  {Con\-temp.\wb Math.}
 \def\comp  {Com\-mun.\wb Math.\wb Phys.}
 \def\cpma  {Com\-pos.\wb Math.}
 \def\duke  {Duke\wB Math.\wb J.}
 \def\fiic  {Fields\wB Institute\wB Commun.}
 \def\foph  {Fortschritte\wB d.\wb Phys.}
 \def\imrn  {Int.\wb Math.\wb Res.\wb Notices}
 \def\inma  {Invent.\wb math.}
 \def\injm  {Int.\wb J.\wb Math.}
 \def\jgap  {J.\wb Geom.\wB and\wB Phys.}
 \def\jktr  {J.\wB Knot\wB Theory\wB and\wB its\wB Ramif.}
 \def\joac  {J.\wB Al\-ge\-bra\-ic\wB Com\-bin.}
 \def\joal  {J.\wB Al\-ge\-bra}
 \def\jomp  {J.\wb Math.\wb Phys.}
 \def\jopa  {J.\wb Phys.\ A}
 \def\jpaa  {J.\wB Pure\wB Appl.\wb Alg.}
 \def\maan  {Math.\wb Annal.}
 \def\nupb  {Nucl.\wb Phys.\ B}
 \def\phep  {Proc.\wb HEP$\!$} 
 \def\phlb  {Phys.\wb Lett.\ B}
 \def\phrl  {Phys.\wb Rev.\wb Lett.}
 \def\ptrs  {Phil.\wb Trans.\wb Roy.\wb Soc.\wB Lon\-don}
 \def\rvmp  {Rev.\wb Math.\wb Phys.}
 \def\sema  {Selecta\wB Mathematica}
 \def\slnm  {Sprin\-ger\wB Lecture\wB Notes\wB in\wB Ma\-the\-matics}
 \def\tams  {Trans.\wb Amer.\wb Math.\wb Soc.}
 \def\topo  {Topology}
 \def\trgr  {Trans\-form.\wB Groups}
 \def\BIR    {{Birk\-h\"au\-ser}}
 \def\OUP    {{Oxford University Press}}
 \def\SV     {{Sprin\-ger Ver\-lag}}
 \def\WS     {{World Scientific}}
 \def\Bo     {{Boston}}
 \def\Si     {{Singapore}}
 \def\NY     {{New York}}

\end{document}